\documentclass[preprint,12pt]{elsarticle}

\usepackage{amsmath}
\usepackage{amssymb} 
\usepackage{bm}
\usepackage{algorithm} 
\usepackage{algorithmic}
\usepackage{listings}   
\usepackage{geometry}
\usepackage{setspace}
\usepackage{titlesec}
\usepackage{titletoc}
\usepackage{graphicx}
\usepackage{indentfirst}
\usepackage{multirow}
\usepackage{enumitem}
\usepackage[htt]{hyphenat}
\usepackage[T1]{fontenc}
\usepackage{booktabs}

\newtheorem{remark}{Remark}
\titleformat{\paragraph}[runin]{\normalfont\normalsize\bfseries}{\theparagraph}{1em}{}
\titlespacing*{\paragraph}{0pt}{3.25ex plus 1ex minus .2ex}{1em}

\newcommand{\xiaosihao}{\fontsize{12pt}{\baselineskip}\selectfont}

 \newcommand{\Rmnum}[1]{\uppercase\expandafter{\romannumeral #1}}

\xiaosihao                                       

\titlecontents{chapter}[0em]{}{\makebox[4.1em][l]{第\xCJKnumber{\thecontentslabel}章}}{}{\titlerule*[0.7pc]{.}\contentspage}

\newcommand{\bfn}{{\boldsymbol n}}

\newcommand{\bfv}{{\boldsymbol v}}

\newcommand{\bfx}{{\boldsymbol x}}

\newcommand{\bfA}{{\boldsymbol A}}

\newcommand{\bfF}{{\boldsymbol F}}

\newcommand{\bfJ}{{\boldsymbol J}}

\newcommand{\mcT}{{\mathcal T}}

\newcommand{\mbR}{{\mathbb R}}

\newcommand{\bgamma}{{\bm \gamma}}

\newcommand{\blambda}{{\bm \lambda}}

\newcommand{\bPhi}{{\bm \Phi}}

\usepackage[scaled=0.9]{FiraMono} \usepackage{comment}
\usepackage{xcolor}

\newcommand{\dd}{\,{\rm d}}

\newcommand{\funding}[1]{\fntext[funding]{#1}}

\usepackage{hyperref}
\definecolor{urlcolor}{HTML}{0f90bf}
\hypersetup{
	plainpages=false,
    colorlinks,
    linktocpage=true,
    linkcolor={blue!90!black},
    citecolor={red!60!black},
    urlcolor={urlcolor}
}

\usepackage{amssymb}
\usepackage{amsmath}

\journal{Neural Networks}

\begin{document}

\begin{frontmatter}

\title{
  Solver-in-the-Loop joint operator learning: fractional Laplace-Beltrami features for interface reconstruction
}

\author[label1]{Yangyang Zheng}
\author[label2]{Huayi Wei}
\ead{weihuayi@xtu.edu.cn}

\author[label3]{Shuhao Cao}
\ead{scao@umkc.edu}

\author[label4]{Ruchi Guo\corref{cor1}}
\ead{ruchiguo@scu.edu.cn}

\cortext[cor1]{Corresponding author}

\affiliation[label1]{organization={School of Mathematics and Computational Science},
            addressline={Xiangtan University},
            city={Xiangtan},
            state={Hunan},
            country={China}}

\affiliation[label2]{organization={School of Mathematics and Computational Science},
            addressline={Xiangtan University, National Center for Applied Mathematics in Hunan, Hunan Key Laboratory for Computation and Simulation in Science and Engineering},
            city={Xiangtan},
            postcode={411105},
            state={Hunan},
            country={China}}

\affiliation[label3]{organization={School of Science and Engineering, University of Missouri-Kansas City},
            city={Kansas City},
            state={Missouri},
            country={United States}}

\affiliation[label4]{organization={School of Mathematics, Sichuan University},
            city={Chengdu},
            state={Sichuan},
            country={China}}

\funding{The first and second authors were supported by the National Natural Science
Foundation of China (NSFC) (Grant Nos. 12371410, 12261131501)
and the construction of innovation provinces in Hunan Province (Grant No. 2021GK1010). The third author was supported in part by NSF grant DMS-2309778.
The fourth author was supported by
National Key R\&D Program of China 2025YFA1018700,
NSFC grant 12571436,
the Fundamental Research Funds for the Central Universities,
the start-up grant of Sichuan University,
NSF grant DMS-2012465.}

\begin{abstract}
In this work, we propose a joint operator learning method for reconstructing images of conductivity coefficients from boundary data.
Inspired by the idea of employing partial differential equation (PDE) solvers as preconditioners for this inverse problem,
we investigate a ``solver-in-the-loop'' training mechanism.
It allows the interaction of learnable parameters integrated in a PDE solver module and those in neural networks for reconstructing images.
Specifically, we employ a fractional Laplace-Beltrami operator with a learnable fractional order, which transforms boundary data into high-dimensional features.
These features then serve as input to a neural network, significantly improving reconstruction accuracy.
For this purpose, a Learning-Automated FEM (LA-FEM) package, facilitating this ``solver-in-the-loop'' property, is developed with PyTorch as a backend.
The new LA-FEM module conveniently allows the auto-differentiation regarding an objective function to freely propagate through the PDE solver from the forward problem and the coupled neural networks for the inverse problem.
\end{abstract}


\end{frontmatter}

\section{Introduction}

Partial differential equations (PDEs) arise in many scientific and engineering
applications.  In PDEs, the usual forward problem concerns solving for the solutions
from the parameters--coefficients, boundary conditions, source and so on.  In
contrast, inverse problems seek to recover those parameters from
extra observed data.  A well-known example is coefficient reconstruction, where
the goal is to infer internal material properties, such as conductivity or
diffusion coefficients, merely from boundary measurements. 
In this work, we focus on Electrical Impedance Tomography (EIT) that is a
specific instance of a coefficient reconstruction problem. EIT aims to
reconstruct the electrical conductivity distribution, usually treated as images, inside a domain from
voltage and current measurements taken on the boundary, respectively. 
It is a non-invasive imaging
technique utilized in various fields, such as medical diagnostics and geophysical exploration. 

In this work, we restrict ourselves to the
scenario of piecewise-constant conductivity values.
This setup is common in many real-world applications of biomedical EIT \cite{2019MartinsSato,2018LiuKhambampatiDu,2021RenWangDong}, 
such as bedside lung monitoring where lung regions and non-lung background can be approximated by uniform conductivity, 
and breast tumor detection where tumor inclusion and healthy background are also likewise well modeled as regions with uniform conductivity.
Moreover, the primary focus of this work is the setting where the conductivity values are (approximately) known a priori, 
while the geometry of the discontinuity is unknown. 
This assumption is also adopted in the aforementioned applications.
More complicated cases of continuous conductivity distribution will be illustrated by numerical examples.
We shall use this problem to demonstrate how PDE solvers can be integrated into a deep learning (DL)
loop and serve as a learnable module for preconditioning and thus to enhance reconstruction.

It is known that obtaining a high-quality reconstruction for inverse problems is
challenging due to their severe ill-posedness. Current algorithms can be categorized into two families. 
The first one is based on iterative algorithms to minimize a data-mismatching function, see e.g., 
\cite{martin1997fem,vauhkonen1999three,jin2017convergent,GuoLinLin2017,2004ChanTai,chung2005electrical,2001ItoKunischLi,2019JinXu,2018GongSchullcke,2020JauhiainenKuusela}.  Alternatively, the second class, called direct methods, aims to build the explicit formula for one-shot reconstruction, such as
factorization methods \cite{kirsch2008factorization}, multiple signal
classification algorithms \cite{ammari2007polarization,
ammari2004reconstruction}, and direct sampling methods (DSMs)
\cite{chow2014direct,2018ChowItoZou,ItoJinZou}. More recently, \cite{2025ItoJinWangEtAlSIIMS} developed a novel iterative DSM by rewriting the formulation of an index function as a variational problem, which combines the strengths of the two.

In addition, DL methods have shown significant
potential for solving inverse problems, and their distinguished feature is the
application of deep neural networks (DNNs) for approximation, replacing
traditional closed-form mathematical formulas.  Supported by data, the DL methods enable
more flexible modeling of complex relationships, resulting in better
reconstructions, especially when empirical structures, unknown hyperparameters,
or restrictive assumptions are largely involved.  Such properties make DL
approaches particularly favorable in inverse problems. For instance,
physics-informed neural networks (PINNs) have been widely used for solving
inverse problems
\cite{2022GaoZahrWang,2022JagtapMaoAdamsKarniadakis,lu2021physics,2019RaissiPerdikarisKarniadakis,2021YangMengKarniadakis,2025WuDuanSunYu}.
We also refer readers to \cite{2024ZhangChanTai} for employing PINN to vessel simulation with complex geometry.
Moreover, the authors in \cite{2020BaoYeZangZhou} proposed a weak adversarial network for solving EIT.
Additionally, as for EIT problems, notable examples of coupling
classical reconstruction methods and DL methods include
\cite{hamilton2019beltrami,hamilton2018deep,2019FanBohorquezYing,fan2020solving,2025WangDengLiu,2025SmylTallmanHomaFlournoy}.  
Moreover, a deepDSM is proposed in \cite{2020GuoJiang,2021JiangLiGuo} that learns local convolutional kernels mimicking the gradient operator of the classical DSMs, and DSM-based DL methods are then widely studied in
\cite{2023LiLiangWang,2023NingHanZou,2025NingHanZou,2025NingZou}.

Interestingly, DL approaches share a conceptual
similarity with traditional methods. One class of DL methods aims to
reconstruct a single instance of the unknown quantity by approximating it with
NNs, a strategy known as ``function learning'' or ``function representation''.  The second class, referred
to as ``operator learning'' (OpL), focuses on approximating the mapping between
the input data and the desired reconstruction with NNs, capturing the underlying
relationship between them.  The approach to approximate a single instance is similar to the classical iterative methods.
In contrast, the OpL approach is close to the classical direct methods, whose evaluation bypasses numerous iterations in the reconstruction process. In classical direct methods, the reconstruction formulas are largely heuristically guided by Green's functions or adjoint formulations, yet sometimes not rigorously derived from a complete inversion of the forward model. OpL replaces these intuitive yet effective mathematical formulas by operator-valued NNs coupled with its training procedure.
The OpL methodology has been widely adopted for solving PDE-related problems 
\cite{2023CaoRoseberryJha,2022DengShinLuZhangKarniadakis,2021LiKovachkiAzizzadenesheli,2023MolinaroYangEngquistMishra,2025YanXuMa}
by learning the inverse or forward operators directly from data.

\begin{figure}[h]
    \centering
    \includegraphics[width=0.7\textwidth]{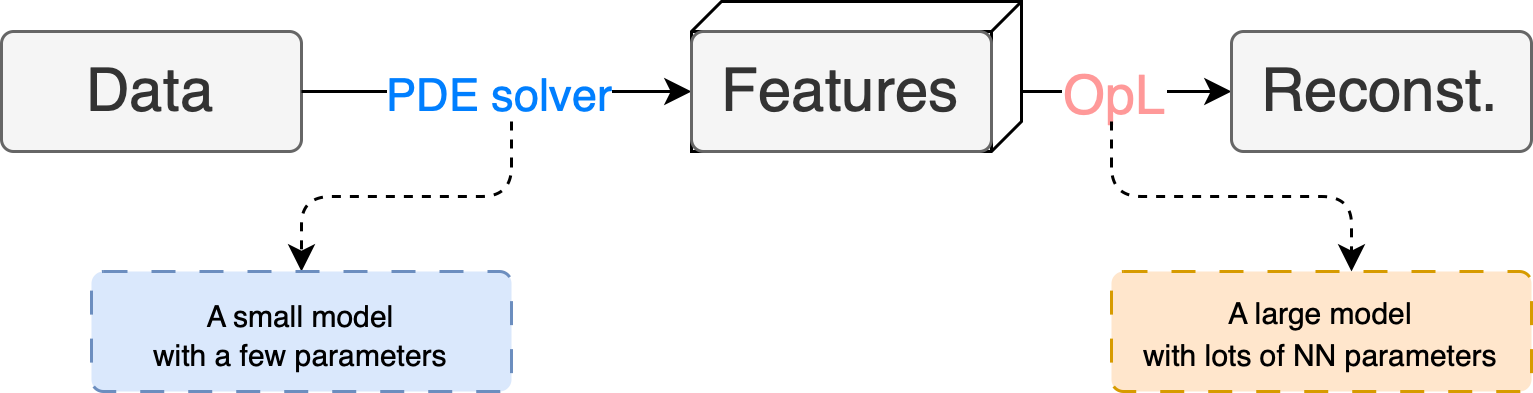}
    \caption{The outline of the $\gamma$-deepDSM.
    The PDE model is more like a ``white box'' containing only a few parameters with very clear mathematical and physical meaning,
    and its main purpose here is to process the data for preconditioning,
    while the OpL model is a ``black box'' containing a large number of network parameters,
    which is to map the processed data to the target images.
    All these parameters are data-driven.}
    \label{fig:outlne}
\end{figure}

For solving ill-posed problems,
preconditioners play a crucial role as they can alleviate ill-posedness to a certain extent \cite{chow2014direct,martin1997fem}.
Despite this, preconditioners are rarely applied in current OpL methods,
due to (i) for complicated nonlinear inverse problems and NNs, development of preconditioners is highly challenging,
and (ii) preconditioners usually incorporate certain PDE solvers which have not been effectively and efficiently integrated with NNs.

In this work, we propose a novel workflow that
incorporates a learnable PDE solver as a module within an OpL architecture.  
The solver contains two key operators:
one is a fractional-order Laplace-Beltrami (fLB) operator with learnable orders ($\gamma$), 
treated as a parameterization of preconditioners for the boundary data,
and the other one is an extension operator of the boundary data to the whole domain.
Conventionally, these fractional orders are hyperparameters empirically
determined and fixed in every scenario in an ad-hoc manner. 
In our approach, they are instead learned optimally in a ``solver-in-the-loop'' manner through the auto-differentiable PDE solver.
The entire solver can be also understood as a data-feature generation mechanism,  
where the feature dimension is much larger than the input data dimension, as shown in Fig. \ref{fig:outlne}.
The features are subsequently mapped to the reconstructed images through NNs.
We name the resulting method $\gamma$-deepDSM.

Recovering fractional orders has been widely studied in inverse problems, see \cite{2021JinZhou,2022JingYamamoto} for example,
and how to leverage them in a data-driven fashion has obtained considerable research attention, 
see \cite{2022GuoWuYuZhou,2019PangLuKarniadakis} for instance. 
However, the methodology here conceptually differs from those in the literature,
in that the fractional orders are treated as learnables for adjusting preconditioners of the OpL models. 
In this sense, it shares a conceptual similarity with those learning regularization parameters, 
see \cite{2021AfkhamChungChung} for instance.

Implementing this strategy requires integrating traditional numerical PDE solvers, such as Finite Element Methods (FEM) or Finite Difference Methods (FDM), with PyTorch. The integration has raised much attention recently.  
For instance, Torch-FEniCS
\cite{2019torchfenics} is an interface to PyTorch, allowing data (solution, source term) to be passed as \texttt{torch.Tensor} to FEniCS modules. 
\cite{2019torchfenics} is implemented taking advantage of FEniCS-adjoint package~\cite{2013FarrellHamFunkeEtAlSISC}, which supports auto-differentiation (AD) for problems like shape optimizations using FEM solvers. However, Torch-FEniCS has to pass numpy arrays to PyTorch \texttt{nn.Module} with a non-native implementation.\footnote{Torch-FEniCS uses the \texttt{torch.from\_numpy} function, which converts numpy's \texttt{ndarray} to \texttt{torch.Tensor}, see Line 69 here: \url{https://github.com/barkm/torch-fenics/blob/master/torch_fenics/torch_fenics.py\#L69}, that cannot preserve the gradient nor the computational graph, if any, associated with the original numpy tensor, unless the gradient tensors are mannually assigned along withtheir graph history in PyTorch, which is vastly inferior to a native implementation.}

Additionally, in \cite{2023femnets,2025LiXiangChenHeaneyDargaville,1994TakeuchiKosugi}, FE shape functions are represented by NNs,
which benefits both applications and theoretical understandings.
A close work is PyTorch-FEA \cite{2023pytorchfea}: a FEM solver
package incorporates PyTorch's
AD for both forward and inverse problems, 
but it still relies on looping operations during matrix assembly.
Magnum.np
\cite{2023magnumnp} implements FDM on PyTorch and also
benefits from the AD to solve inverse problems. 
In \cite{2025NeuralSEM}, a PINN is combined with a high-order spectral element method to obtain a robust simulation of
Navier-Stokes equations.
Notably, MFEM \cite{2021mfem} enables flexible PDE
solutions through customizable components like mesh, discretization,
solvers and AD, supported by GPUs.
However, to the best of our knowledge, MFEM has not been upgraded to incorporate DL models so far.

In this work, based on FEALPy \cite{wei2024fealpy}, 
we take an alternative approach compared to frameworks like Torch-FEniCS or PyTorch-FEA by implementing a modular design for FEM computations, natively integrated with PyTorch's tensor computing and AD features.
Moreover, it utilizes high-level interfaces while with low-level access to
manipulate data.
FEALPy is a cross-platform Computed-Aided-Engineering engine, 
supporting multiple computation backends including PyTorch.
Its design structure already facilitates many applications, see \cite{2022ChenHuangWei,2018WangWei} for example.
The novel modules in FEALPy, named Learning-Automated FEM (LA-FEM), 
offer several remarkable features:
\begin{itemize}[leftmargin=2em] 
    \item All LA-FEM modules are compatible with PyTorch's AD interface, allowing all FEM components to incorporate learnable parameters. Optimizing
        these parameters mirrors the training process of NNs. 
    \item LA-FEM is fully tensor-based, encompassing FE shape function
        generation, integration, and matrix assembly, and harnesses the parallel
        computing power of GPUs.  
    \item LA-FEM supports native batch-wise computation, consistent with
        PyTorch's data manipulation, enabling simultaneous solutions of multiple
        PDEs in a parallelized framework.  
\end{itemize}
{In the present work, all these features collectively aid in realizing the desired ``solver-in-the-loop" method for training the fLB orders to produce an optimal preconditioner for improving the reconstructions.}

Moreover, LA-FEM is structured around highly object-oriented modules, including
\texttt{Mesh}, \texttt{FunctionSpace}, and \texttt{FEM}, with the \texttt{FEM}
module further composed of \texttt{Integrator} and \texttt{Form} objects. This
modular architecture makes LA-FEM highly scalable, maintainable, and
user-friendly, supporting customized FEM computations across various
applications.

In the next section, we present the architecture of $\gamma$-deepDSM. 
In Section \ref{sec:LAFEM}, we discuss the implementation details of LA-FEM.
In Section \ref{sec:num}, we present the numerical results.
In the last section, we draw the conclusion.
Usage samples and architecture details of LAFEM are given in Appendix.

\section{The \texorpdfstring{$\gamma$}{Gamma}-deepDSM}

Consider a bounded domain $\Omega \subset \mbR^2$ and an inclusion $\Omega_1\subseteq\Omega$ filled with different conductive media.
The conductivity of the background homogeneous medium $\Omega_0 := \Omega\backslash \Omega_1$ and $\Omega_1$ are denoted by $\sigma_0$ and $\sigma_1$, respectively.
Let $\sigma$ be a piecewise-constant function such that $\sigma|_{\Omega_1}=\sigma_1$ and $\sigma|_{\Omega_0}=\sigma_0$.
The following $L$ forward model equations are obtained by injecting $L$ different electrical currents to the boundary:
\begin{equation}
\begin{aligned}
-\nabla\cdot(\sigma\nabla u_l)=0, \quad \text{in }\Omega,
~~~~~~
\text{with}
~~
\sigma\frac{\partial u_l}{\partial n}=g_{N,l}, \quad \text{on }\partial\Omega, \quad l=1,2,\dots,L,
\end{aligned} \label{eq_eit}
\end{equation}
where $g_{N,l} \in H^{-1/2}(\partial\Omega)$ is the current density such that
$
\int_{\partial \Omega}g_{N,l}(s)\text{d}s=0.
$
For EIT, $\sigma$ is unknown and needs to be recovered by the Neumann and Dirichlet boundary data,  
known as Cauchy data pairs $\{(g_{N,l}, g_{D,l})\}_{l=1}^L$ with $g_{D,l}:=u_l|_{\partial\Omega}$ being the voltage.
Given a $\sigma$, the Neumann-to-Dirichlet mapping $\Lambda_{\sigma}$ from $g_{N,l}$ to $g_{D,l}$ (NtD mapping) can be uniquely determined by solving \eqref{eq_eit} and applying the trace operator. 
Given the full ``knowledge'' of $\Lambda_{\sigma}$ by knowing how $\{g_{D,l}\}_{l=1}^{\infty}$ respond to only countably many $\{g_{N,l} \}_{l=1}^{\infty}$ (since $H^{-1/2}(\partial\Omega)$ is Hilbert), 
the unknown $\sigma$ can be recovered theoretically \cite{1988Alessandrini}.
However, this procedure is very ill-posed, only log-stable with respect to perturbation according to \cite{1988Alessandrini}.
It becomes even more ill-posed since merely a finitely many $L$ Cauchy pairs~\cite{2016CaroFerreiraRuizJDE} are known in this paper, $\Lambda_{\sigma}$, i.e., $\{(g_{N,l}, g_{D,l})\}_{l=1}^L$. In our simplified setting (piecewise constant conductivity), it can be shown that the inverse problem is Lipschitz stable~\cite{2022AlbertiSantacesariaAfRMA,2005AlessandriniVessellaAAM}.

In the following discussion, we define the boundary data for a series of scattered potentials as \begin{equation}
\label{bc_data_eq1}
\xi_l =g_{D,l}-\Lambda_{\sigma_0}(g_{N,l}) = (\Lambda_{\sigma}-\Lambda_{\sigma_0}) g_{N,l}, ~~~~ l=1,...,L.
\end{equation}
It characterizes the discrepancies between the voltages with inclusions and those without. In the next subsection, we will establish the nonlinear equation illustrating how these scattered data encode the conductivity $\sigma$.
For simplicity, we denote $\Xi = \{ \xi_l \}_{l=1}^L$.

\subsection{The nonlinear equation of EIT}
Now, we follow the procedure recently introduced in \cite{2025ItoJinWangEtAlSIIMS} to formulate an equation of $\sigma$.
To this end, we only consider the case of a single Cauchy data pair $(g_N, g_D)$.
To avoid ambiguity with $\sigma$ reconstruction using full measurement, we introduce $\sigma_{\dag}$ and $u_{\dag} = u(\sigma_{\dag})$ as the true coefficient and the corresponding solution, respectively, to a data pair in \eqref{eq_eit}.
 Denote $(\cdot,\cdot)_{\Omega}$ as the $L^2(\Omega)$-inner product; 
 then, it is not hard to see that the weak form of \eqref{eq_eit} can be equivalently written as
\begin{equation}
\label{eq_eit2}
\begin{aligned}
( \sigma_0\nabla u_{\dag}, \nabla v )_{\Omega} + ( \tilde{\sigma}_{\dag}\nabla u_{\dag}, \nabla v )_{\Omega} = ( g_{N}, v )_{\partial \Omega}, ~~~~~ \forall v\in H^1(\Omega),
\end{aligned}
\end{equation}
with $\tilde{\sigma}_{\dag} = \sigma_{\dag} - \sigma_0$. 
In the following discussion, we denote $\langle\cdot,\cdot\rangle$ as the standard $H^{-1}$--$H^1$ duality product
and define $\overline{H}^1(\Omega) = \{ u\in H^1(\Omega): (u,1)_{\partial\Omega} =0 \}$.
Introduce the following linear operators:
\begin{equation}
\label{eq_eit3}
\begin{aligned}
&\mathcal{A} ~ : ~ \overline{H}^1(\Omega) \rightarrow H^{-1}(\Omega)
~~~
\text{with}
~~~
\langle \mathcal{A}u, v \rangle := ( \sigma_0\nabla u, \nabla v )_{\Omega}, \\
&\mathcal{B}[\tilde{\sigma}]~ : ~ H^1(\Omega) \rightarrow H^{-1}(\Omega)
~~~
\text{with}
~~~
\langle \mathcal{B}[\tilde{\sigma}]u, v \rangle := ( \tilde{\sigma}\nabla u, \nabla v )_{\Omega},
\end{aligned}
\end{equation}
and $g_N'$ as certain functional-valued extension of $g_N$
\begin{equation}
\label{eq_gN_ext}
g'_N \in H^{-1}(\Omega) 
~~~
\text{with}
~~~ 
\langle g'_N, v \rangle :=  ( g_{N}, v )_{\partial \Omega}.
\end{equation}
Additionally, we will use the following adjoint operators:
\begin{equation*}
( \mathcal{B}_{\tau}[u] \tilde{\sigma}, v )_{\Omega} =  ( \tilde{\sigma},   \mathcal{B}^*_{\tau}[u]v )_{\Omega}  = ( \tilde{\sigma}\nabla u, \nabla v )_{\Omega},
 ~~~ \forall u,v,\tilde{\sigma}.
\end{equation*}
Then, it is obvious that
\begin{equation}
\label{eq_eit5}
 \mathcal{B}^*_{\tau}[u]v = \nabla u \cdot\nabla v.
\end{equation}
Now, \eqref{eq_eit2} can be reformulated as an operator equation:
\begin{equation*}
u_{\dag} = \mathcal{A}^{-1}\left( g'_N -   \mathcal{B}_{\tau}[u_{\dag}] \tilde{\sigma}_{\dag}\right).
\end{equation*}
Introducing the trace operator $\mathcal{T}: H^{1}(\Omega)\rightarrow H^{1/2}(\partial \Omega)$,
we obtain
\begin{equation*}
g_{D} - \mathcal{T} \mathcal{A}^{-1} g'_N = - \mathcal{T} \mathcal{A}^{-1}  \mathcal{B}_{\tau}[u_{\dag}] \tilde{\sigma}_{\dag} .
\end{equation*}
As $\Lambda_{\sigma_0} g_N = \mathcal{T} \mathcal{A}^{-1} g'_N$ by \eqref{eq_gN_ext},
the left-hand side above is exactly the scattered data $\xi := g_{D} -  \Lambda_{\sigma_0} g_N $,
which yields a nonlinear operator equation for $\sigma_{\dag}$:
\begin{equation}
\label{eq_eit7}
\xi = - \mathcal{F}[\sigma_{\dag}] \sigma_{\dag} ~~~~~~~ \text{with} ~~~ 
\mathcal{F}[\sigma_{\dag}] := \mathcal{T} \mathcal{A}^{-1}  \mathcal{B}_{\tau}[u(\sigma_{\dag})]. \end{equation}
When thinking in a discrete sense, $\mathcal{F}$ maps high-dimensional interior data to lower-dimensional boundary data, solving \eqref{eq_eit7} for $\sigma_{\dag}$ is severely ill-posed. 
With a slight abuse of the presentation order, we introduce an operator $\mathfrak{L}$ of which its meaning and specific form will be specified later. Denoting $\mathcal{R} := \mathcal{F}^{*}\mathfrak{L}\mathcal{F}$,
Using the symmetry of $\mathcal{A}$, we deduce from \eqref{eq_eit7} that
\begin{equation}
\label{eq_approx_1}
\mathcal{R} \sigma_{\dagger} = -\mathcal{F}^* \mathfrak{L} \xi  = -\mathcal{B}^*_{\tau}[u_{\dag}] \mathcal{A}^{-1} \mathcal{T}^* \mathfrak{L} \xi,
\end{equation}
where $\mathcal{T}^*: H^{-1/2}(\partial \Omega) \to H^{-1}(\Omega)$ is the adjoint of the trace operator.
Denote $\phi = \mathcal{A}^{-1} \mathcal{T}^* \mathfrak{L} \xi$, it is straighforward to see that by definitions above $\phi$ satisfies
\begin{equation}
\label{eq_approx_2}
-\nabla\cdot(\sigma_0 \nabla \phi) = 0, ~~~~ \nabla \phi\cdot \bfn =  \mathfrak{L} \xi.
\end{equation}
Combining \eqref{eq_eit5}, \eqref{eq_approx_1} and \eqref{eq_approx_2}, we arrive at the identity
\begin{equation}
\label{eq_approx_3}
\mathcal{R} \sigma_{\dagger}  = -\nabla u_{\dag}\cdot\nabla \phi,
~~~~ \text{with} ~~
\mathcal{R} = \mathcal{B}^*_{\tau}[u_{\dag}] \mathcal{A}^{-1} \mathcal{T}^*  \mathfrak{L} \mathcal{T} \mathcal{A}^{-1}  \mathcal{B}_{\tau}[u_{\dag}].
\end{equation}
Notice that, in actual computation involving discretizations using this reformulation, $\mathcal{R}$ may not be invertible in general, and this nature is inherited from the severe ill-posedness from the original formulation of EIT (from the NtD map to the conductivity). 
In the classical DSM such as \cite{chow2014direct}, the index function $\mathcal{I}(\cdot)$ can be viewed as certain inversion to an empirically constructed $\widetilde{\mathcal{R}} \approx \mathcal{R}$
\begin{equation}
\label{eq_approx_5}
\mathcal{I}(x) : = \widetilde{\mathcal{R}}^{-1} |\nabla \phi(x) |^2 \approx \sigma_{\dag},
\end{equation}
where $\nabla \phi$ is constructed by a harmonic extension to approximate $\nabla u_{\dag}$. DSM is extremely efficient in that, even for a single Cauchy data pair, one only has to sample a few points to get a good estimate of the location of the inclusion. 

However, a simple formula, such as \eqref{eq_approx_5}, still suffers from limited accuracy (the reconstruction formula uses a single term in an asymptotic expansion in \cite{cedio1998identification}), and it is unclear how to generalize the formula to multiple data pairs, as well as how to choose a suitable preconditioner $\mathfrak{L}$.
The present work aims to address these issues by integrating a learnable FEM solver and an NN, 
where the former one is used to learn the preconditioner, while the latter one is to approximate nonlinear inverse operator in a data-driven manner.

One major component of $\mathcal{R}$ is the linear operator $\mathcal{A}^{-1/2}\mathcal{T}^*\mathfrak{L} \mathcal{T} \mathcal{A}^{-1/2} = (\mathcal{A}^{-1/2}\mathcal{T}^*\mathfrak{L}) (\mathcal{T} \mathcal{A}^{-1/2})$, making it spectrally equivalent to $(\mathcal{T} \mathcal{A}^{-1/2})(\mathcal{A}^{-1/2}\mathcal{T}^*\mathfrak{L}) = \mathcal{T}\mathcal{A}^{-1}\mathcal{T}^* \mathfrak{L}$ due to the symmetry of $\mathcal{A}$.
Notice that $\mathcal{T}\mathcal{A}^{-1}\mathcal{T}^* = \Lambda_{\sigma_0}$ in the sense of \eqref{eq_approx_2}. Moreover, in certain simple cases (e.g., on the boundary of the 2D unit disk with only $\sigma_0$ as the conductivity in the interior, see e.g., \cite{1995PidcockKuzuogluLeblebicioglu}), it behaves like the fractional Laplace-Beltrami (fLB) operator $(-\Delta_{\partial\Omega})^{-1/2} = d/d\theta$.
In this regard, $\mathfrak{L}$ can be chosen to serve as a preconditioner for this linear operator $\mathcal{A}^{-1/2}\mathcal{T}^*\mathfrak{L} \mathcal{T} \mathcal{A}^{-1/2}$ in helping inverting $\mathcal{R}$, and one typical choice is its inverse $(-\Delta_{\partial\Omega})^{1/2}$ as in the original DSM paper~\cite{chow2014direct}. In general, we may use $(-\Delta_{\partial\Omega})^{-\gamma}$ for some $\gamma>0$, despite that the optimal order $\gamma$ is unknown in general (see Section \ref{subsec:fLBanalysis}).
In this work, we shall use fLB with learnable order to parameterize $\mathfrak{L}$ for improving the reconstruction accuracy.

\subsection{The fractional Laplace-Beltrami operator}

In this subsection, we first introduce a regularized fLB operator and describe its preconditioning property.
For a function $f: M\to \mathbb{R}$ on a general closed Riemann manifold $M$, the LB operator is defined as
$$
\Delta_{M}: f \mapsto \Delta_{M} f:= \nabla_M\cdot(\nabla_M f),
$$
where $\nabla_M$ represents the gradient operator on the manifold $M$.
Among various definitions of fractional derivatives~\cite{2016BucurValdinoci,2021Jin,2017Kwanicki,2020LischkePangGulianSong}, 
we opt for the Fourier-transform definition in our method,  
as the LB operator is a self-adjoint operator on a closed surface $\partial\Omega$ with no boundary condition imposed. 
Since $-\Delta_{\partial \Omega}$ is semi-positive-definite, for which zero is a simple eigenvalue associated with a constant eigenfunction,
we let $\lambda_k, \psi_k$ be its (non-zero) eigenvalues and orthonormal eigenfunctions, respectively.
Then, every function $u\in H^1(\partial \Omega)$ admits the following expression 
$$
u = \sum_{k=1}^\infty \alpha_k \psi_k,
$$
and as $-\Delta_{\partial \Omega}$ is bounded from $H^1$ to $H^{-1}$, we have
$$
-\Delta_{\partial \Omega}u 
= -\Delta_{\partial \Omega}\left(\sum_{k=0}^\infty \alpha_k \psi_k \right) 
= \sum_{k=0}^\infty \alpha_k \lambda_k \psi_k.
$$
The fLB operator can be then defined as
\begin{equation}
    (-\Delta_{\partial \Omega})^{\gamma} u 
    := \sum_{k=0}^\infty \alpha_k \lambda_k^\gamma \psi_k. \label{frac-lap-c}
\end{equation}
To represent the fLB operator in our method, eigenpairs can be calculated a priori in an offline stage, e.g., one may simply choose to use FEM on a uniform mesh of $\partial\Omega$ to compute the eigenpairs.

\begin{figure}[h]
    \centering
    \includegraphics[width=1.0\textwidth]{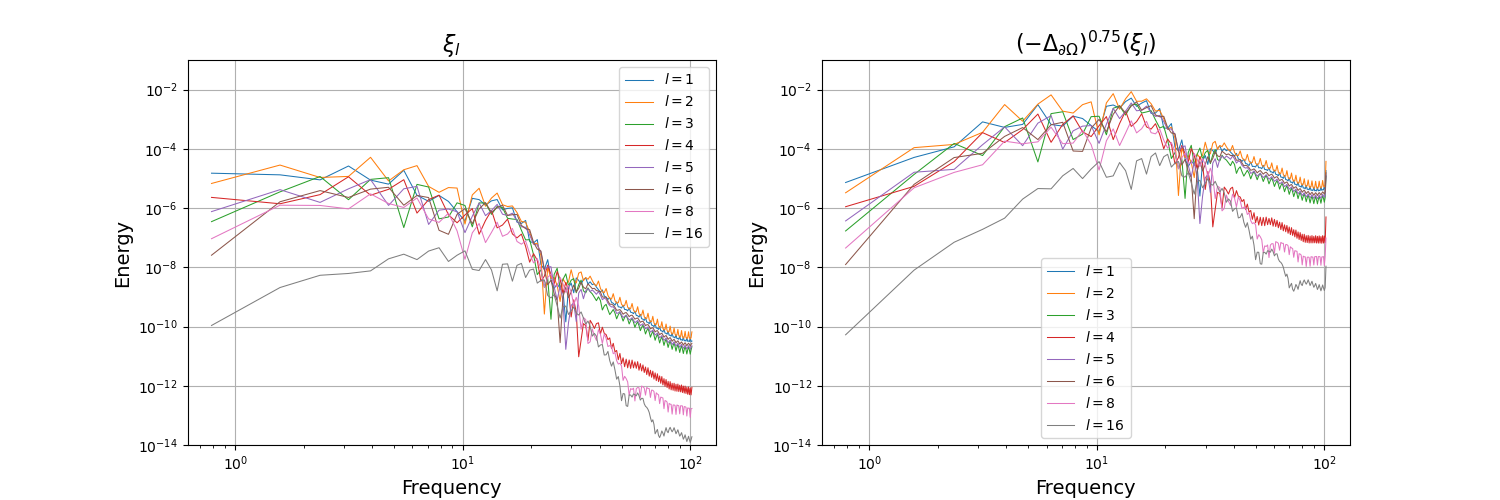}
    \caption{Left: spectrum of channels in $\xi_l$ and $ \mathfrak{L}^{\gamma}_{\partial\Omega} \xi_l$, 
    where $g_{D,l}=\cos(l\theta(\bfx))$, $\bfx\in\partial\Omega$, $l=1,...,6,8$ and $16$.
    Right: the fLB operator with $\gamma=0.75$ effectively increases the energy of high-frequency information.}
    \label{spectrum}
\end{figure}

For $\partial \Omega$ of a simply-connected domain $\Omega$, 
the eigenvalues generally admit the growth order $\lambda_{k} = \mathcal{O}(k^2)$. It is well-known that a typical FEM solver cannot capture the accuracy of eigenpairs associated with very high frequency due to the mesh size and the polynomial order. Thus, with a truncation $K_0\in \mathbb{Z}^+$, we define 
\begin{equation}
\label{frac-lap-c_trunc}
\mathfrak{L}^{\gamma}_{\partial\Omega} u = \sum_{k=0}^{K_0} \alpha_k \lambda_k^{\gamma} \psi_k. 
\end{equation}
For $u\in H^1(\partial \Omega)$, by Plancherel theorem, $\sum_{k=0}^\infty \lambda_k^2 |\alpha_k|^2 < \infty$, 
the energy $E_k=|\alpha|^2_k \sim o(k^{-3-\epsilon})$ in each frequency mode decays very fast as $k$ gets larger. 
Consequently, the high-frequency components carry significantly less energy compared to the low-frequency components, 
making it difficult to extract the information contained in these high-frequency as the SNR decreases. 
However, as $\lambda_k \rightarrow \infty$, we see that, with an extra multiplier $\lambda_k^{\gamma}$, $\lambda_k^{\gamma}\alpha_k$ can mitigate the decay of $\alpha_k$ to a certain extent, e.g., see Fig. \ref{spectrum}.
In this regard, from the perspective of spectrum, the fLB operator can precondition the inverse problem, and thus enhance the reconstruction accuracy.

However, $\gamma$ cannot be too large as $\lambda^{\gamma}_k$ can also magnify the noise of $u$'s measurements. Numerical experiments suggest that the model becomes very sensitive to noise with an overly large $\gamma$, 
significantly deteriorating prediction accuracy for noisy data. 
Therefore, a suitable $\gamma$ is important for both improving the accuracy and enhancing the stability.
In practice, the optimal choice of $\gamma$ depends on several factors, 
including the characteristics of the data, as well as the type and magnitudes of noise involved. In the original DSM paper~\cite{chow2014direct}, the authors gave a heuristic argument and numerically illustrated that opting for a larger $\gamma$ improves the performance of the index function to locate the inhomogeneities. In the following section, we provide a more rigorous exposition in a simplified setting to justify the approach of looking for a $\gamma$ through learning-based approaches.

\subsection{Analysis of the fLB operator on the index functional}
\label{subsec:fLBanalysis}
In this subsection, we shall analyze how the fLB operator can aid in improving a signal-to-noise ratio (SNR) ratio,
by following \cite[Section 4.2]{chow2014direct} to take the special case that $\Omega = B_1(0)$ is the unit disk and one has a single Cauchy data pair $\{g_N, g_D\}$ to form $\xi$ in \eqref{bc_data_eq1}.
With a slight abuse of notation, a function defined on the boundary will be written as a function on $\mathbb{S}^1$, i.e., $\xi = \xi(\theta)$ and $\xi(0) = \xi(2\pi)$.

It is not hard to see that, with the definition in \eqref{frac-lap-c} and \eqref{frac-lap-c_trunc}, the fLB operator would amplify both signal and noise. 
Let us make this more precise.
Consider the measured scattered potential:
\begin{equation}
    \xi^{\delta}(\theta) := \xi(\theta) + \eta(\theta),
\end{equation}
where, with Fourier transformation, the true data $\xi(\theta):= \sum_{k\neq 0} \xi_k e^{ik \theta}$ with $\xi_k$ admitting a spectral decay (e.g., see \cite[Section 4.1--4.2]{chow2014direct}), and the noise $\eta(\theta):= \sum_{k\neq 0} \eta_k e^{ik \theta}$ with $\{\eta_k\}$ i.i.d. Gaussian with $\mathbb{E}|\eta_k|^2 = \delta$.
As both $\xi$ and $\eta$ are real data, there hold $\xi_{-k} = \overline{\xi_k}$ and $\eta_{-k} = \overline{\eta_k}$. 
Using the same modes truncation $K_0$ in \eqref{frac-lap-c_trunc} for $\xi^{\delta}$ (we keep the notation $\xi^{\delta}$ for simplicity):
\begin{equation}
    \xi^{\delta} = \sum_{1\leq |k|\leq K_0} (\xi_k + \eta_k) e^{ik \theta}, \quad \text{ and } \quad
\mathfrak{L}^{\gamma}_{\partial\Omega}  \xi^{\delta} =
 \sum_{1\leq |k|\leq K_0} |k|^{2\gamma} (\xi_k + \eta_k) e^{ik \theta}.
\end{equation}
Then, by the DSM in \cite{chow2014direct}: denote $\phi^{\delta}_\gamma$ as the solution to the following problem
\begin{equation}
\Delta\phi^{\delta}_\gamma=0 \quad\text{in}~\Omega, \quad
\nabla\phi^{\delta}_\gamma\cdot\bfn=\mathfrak{L}^{\gamma}_{\partial\Omega}  \xi^{\delta}\quad\text{on}~\partial\Omega, \quad
\int_{\partial\Omega}\phi^{\delta}_\gamma \dd S=0. \label{phi_formula}
\end{equation}
Using the standard separation of variable approach for Poisson problem on a disk, we have
\begin{equation}
\label{phi_gamma_delta}
\phi^{\delta}_\gamma(r,\theta) = \sum_{1\leq |k|\leq K_0} |k|^{2\gamma-1}(\xi_k + \eta_k) r^{|k|}e^{ik \theta} := \phi_{\gamma} + \phi^{\text{noise}}_\gamma.
\end{equation}

Based on \eqref{phi_gamma_delta}, we define the signal and noise energy:
\begin{equation}
\label{SNenergy}
S(\gamma) = \|\nabla \phi_{\gamma}\|_{L^2(\Omega)}^2, 
~~~~
N(\gamma) = \mathbb{E}\left\|\nabla \phi_{\gamma}^{\text {noise }}\right\|_{L^2(\Omega)}^2.
\end{equation}
With these quantities, we introduce a special ``Signal-to-Noise'' ratio (SNR):
\begin{equation}
\label{def_snr_global}
    \text{SNR} := \frac{ S(\gamma) }{ N(\gamma) } = 
\frac{\sum_{|k|=1}^{K_0} k^{4 \gamma-1}|\xi_{k}|^2}{\delta \sum_{|k|=1}^{K_0} k^{4 \gamma-1}}.
\end{equation}

To analyze this SNR, especially its dependence on $\gamma$, we consider a simple inclusion where $\Omega_1 = B_\rho(0)$, i.e., the inclusion is a smaller disk of radius $\rho \in (0,1)$ centered at the origin inside $\Omega=B_1(0)$, together with that the input current has a single mode, 
say $g_m = \cos{(m\theta)}$ for example, 
$m\in \mathbb{Z}^+$ and $1 \le m < K_0$.
We shall use this case to illustrate how $\text{SNR}_\gamma$ behaves with respect to $\gamma$.
Using the result of the Poincar\'{e}-Steklov eigenvalue problem \cite{1995PidcockKuzuogluLeblebicioglu}, we have
\begin{equation}
   \xi =  (\Lambda_{\sigma} - \Lambda_{\sigma_0}) g_N = (\lambda_m - m^{-1}) \cos{(m\theta)},
\end{equation}
where $\lambda_m$ (with duplication) is an eigenvalue of $\Lambda_\sigma$:
\begin{equation}
\lambda_{ m} = \frac{1}{m}\frac{1 - \rho^{2m}\mu}{1+\rho^{2m}\mu}, \quad \text{ with } \quad \mu= \frac{1-\sqrt{\sigma_1}}{1+\sqrt{\sigma_1}}.
\end{equation}
Namely, $\xi$ only contains $\xi_{\pm m} = \pm (\lambda_m - m^{-1})/2$.
Then, the SNR defined in \eqref{def_snr_global} can be expressed as 
\begin{equation}
\text{SNR}(\gamma,m) = \delta^{-1} C(\rho,\mu,m) \frac{ m^{4 \gamma-1}}{  \sum_{k=1}^{K_0} k^{4 \gamma-1}},
~~~
\text{where} ~
C(\rho,\mu,m) = \frac{  \rho^{4m} \mu^2}{m^2(1+\rho^{2m}\mu)^2}.
\end{equation}
Notice that the function $C$ is independent of $\gamma$.
Define the function $\chi(s) = m^s/ (\sum_{k=1}^{K_0} k^{s})$ and the constant $\overline{K}_0 = (K_0!)^{1/K_0}\approx K_0/e$.
By elementary calculus, one can show that, for $1\le m \le \overline{K}_0$, $\chi(s)$ decreases for $s\in(-1,+\infty)$,
and for $\overline{K}_0 < m \le {K}_0$, $\chi(s)$ first increases and then decreases over $(-1,+\infty)$.
Thus, for each $m$, there uniquely exists one $\gamma^{\star}\in(0,\infty)$ to maximize $\text{SNR}(\gamma,m)$.
However, as $m \rightarrow K_0$, the optimal $\gamma^{\star}(m) \rightarrow \infty$, which may enhance the noise too much.
To model the effect of $\gamma$ on the reconstruction accuracy, we further propose a very simple model by incorporating a ``noise budget'' when maximizing $\text{SNR}$:
\begin{equation}
\label{SNR_noise} 
\max_{\gamma} ~ \text{SNR}(\gamma,m), ~~~~ \text{subject to} ~~ N_m(\gamma) \le N_{\max},
\end{equation}
where $N_m(\gamma)= \delta m^{4\gamma-1}$ is the noise at the frequency mode $m$.
Due to the monotonic increasing property of $N_m$, it is easy to see that the constraint becomes 
\begin{equation}
\label{SNR_noise_2}
\gamma \le \bar{\gamma}(m):= \frac{1}{4} + \frac{\ln(N_{\max}/\delta)}{4\ln(m)}.
\end{equation}
For large enough $m$, this bound $\bar{\gamma}(m)$ would be smaller than $\gamma^{\star}(m)$, and the previous discussion shows that $\text{SNR}(\gamma,m)$ is increasing over $(0,\bar{\gamma}(m))$.
Therefore, the optimizer to the constrained optimization problem \eqref{SNR_noise} is $\bar{\gamma}(m)$,
and it can be understood as an adaptive spectral filter that balances signal and noise across different frequency channels.
This quantity is actually decreasing with respect to $m$,
which agrees with the numerical results observed in Fig. \ref{evo_s}.

Nevertheless, the analytical expression of the optimal $\gamma$ remains largely unclear in general, 
particularly when the inclusions exhibit complex geometries and when multiple Cauchy data pairs are incorporated as different channels in a DNN-parameterized index functional~\cite{2020GuoJiang,2023GuoCaoChen}.

\begin{figure}[h]
    \centering
    \includegraphics[width=\textwidth]{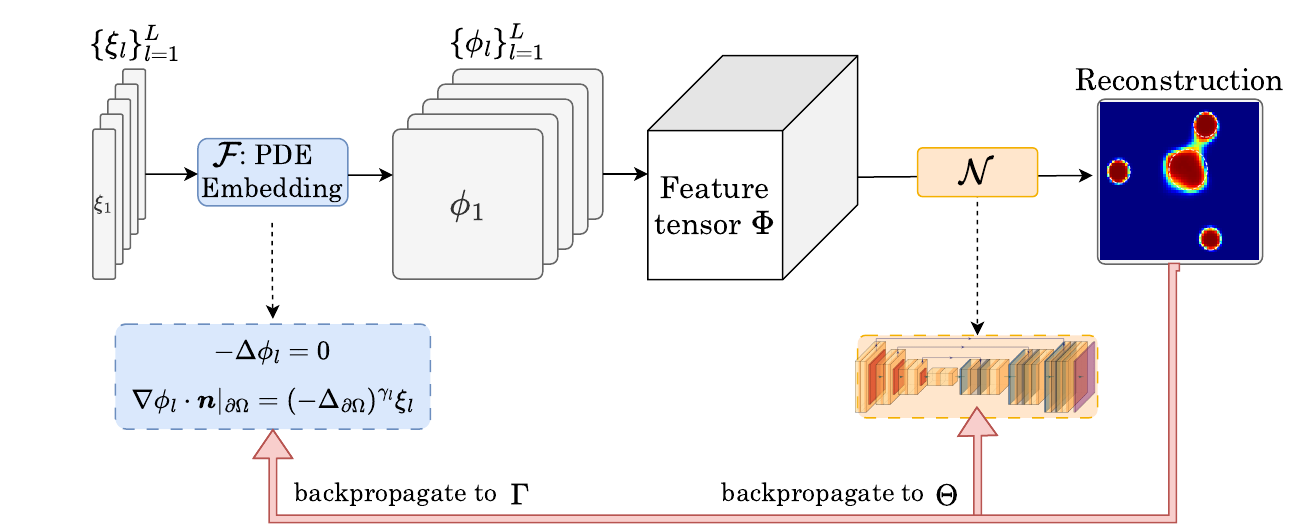}
    \caption{The schematics of $\gamma$-deepDSM: 
    in the first module the learnables are the fractional orders,
    while in the second module the learnables are the network parameters.}
    \label{fig:workflow}
\end{figure}

\subsection{Solver-in-the-loop joint OpL}

We now need to establish a DL workflow from $\Xi=\{\xi_l\}$ to $\sigma$ for reconstruction.
Following \eqref{eq_approx_2}, with the processed data, we define the data features $\phi_l$:
\begin{equation}
\label{bc_data_eq2}
-\Delta\phi_l=0 \quad\text{in}~\Omega, \quad
\nabla\phi_l\cdot\bfn= \mathfrak{L}^{\gamma_l}_{\partial\Omega} \xi_l \quad\text{on}~\partial\Omega, \quad
\int_{\partial\Omega}\phi_l \dd S=0.
\end{equation}
For a given $\gamma_l$, we can obtain $\phi_l$ by a fast Poisson solver using FEM.
Then, these data features $\Phi:=\{ \phi_l \}_{l=1}^L$ are then concatenated into a tensor as the input of an NN, which will be trained to reconstruct the conductivity. 
The highlight of the whole workflow is then to implement an FEM solver with parameters that allow their gradient obtained by AD to be propagated all the way through from the output of the NN. It enables such two components to be trained jointly:
\begin{equation}
\label{workflow_1}
\mathcal{I} = \mathcal{N}(\Phi;\Theta), 
~~~~~~ \Phi = \mathcal{F}( \Xi;\Gamma),
~~~~~~ \mathcal{G}:= \mathcal{N}\circ \mathcal{F},
\end{equation}
where $\mathcal{N}$ is a network with the set of learnable parameters $\Theta$, 
and $\mathcal{F}$ denotes a FEM solver for elliptic equations in \eqref{bc_data_eq2} with $\Gamma = \{ \gamma_1,\gamma_2,...,\gamma_L\}$ being the set of learnable parameters. 
This joint-training method shall be referred to as \textbf{$\gamma$-deepDSM}, and we refer readers to Fig. \ref{fig:workflow} for a schematic diagram.
Assuming that $\sigma$ follows a given distribution $\nu$, 
we need to solve the minimization problem to determine $\mathcal{G}$:
\begin{equation}
\label{loss_fun}
\min_{\Theta,\Gamma} \mathcal{L} (\Theta,\Gamma) := \mathbb{E}_{\sigma \sim \nu} \left[ \mathcal{M}(\mathcal{G}(\Xi;\Gamma,\Theta), \sigma) \right],
\end{equation}
where $\mathcal{M}$ is a certain distance between the NN output $\mathcal{G}(\Xi)$ and $\sigma$.
For the NN module, we employ a standard U-Net~\cite{2015RonnebergerFischerBrox} to fit the nonlinear relationship between the data features $\Phi$ and the inclusion distribution. This process of joint-forward-inverse-operator-learning largely resembles the ``solver-in-the-loop'' approach, e.g., see \cite{2020UmBrandFeiEtAlNIPS,2021KochkovSmithAlievaEtAlPNAS}. 
Nevertheless, note that in these work, only the forward PDE solvers have trainable components,
and these differentiable solvers are not connected with operator learners.
Thus, one major contribution of the present work is to show that learnable PDE solvers can benefit OpL.

\begin{remark}
~\\
\begin{itemize}
\item \textbf{Dimension Lifting}. The operator $\mathcal{F}$ maps the lower-dimensional input boundary data, such as $g_N$, $g_D$ and $\xi$, 
to functions on the whole domain $\Omega$.
This procedure is quite different from many existing encoders that aim to represent spatially higher-dimensional data on a lower-dimensional manifold, e.g., \cite{2021PakravanMistaniAragonCalvoEtAlJCP}. 
A typical example is U-net, which employs learned lower-dimensional latent representation of input data.
The proposed scheme in our method instead seeks for higher-dimensional latent spaces, similar to the ``inverse bottleneck'' technique discussed in \cite{2022LiuMaoWuFeichtenhofe}. In this regard, the module of $\mathcal{F}$ can be interpreted as a special encoder that lifts, instead of reducing, data dimension.
The dimension-lifted images $\Phi$ can be considered as PDE-features with clear mathematical meaning.

\item \textbf{Regularity enhancement}. A useful way to justify this ``lifting'' step is to view it as a geometry-consistent feature map: 
while the harmonic extension does not add new information, 
it re-encodes the boundary data into an interior field whose geometry is better aligned with the target image,
resulting in a structure more consistent with the true operator \cite{chow2014direct} and making the downstream network easier to train.
Moreover, such high-dimensional features $\Phi$ are generated by solving elliptic PDEs, improving training and learning through a regularity enhancement mechanism.
More specifically, boundary data contaminated by white noise only admits Sobolev regularity $H^{s}(\partial \Omega)$ with $s<-(d-1)/2$, $d=2,3$.
Applying the fLB operator may further reduce the regularity to $H^{s-2\gamma}(\partial \Omega)$, 
thereby amplifying high-frequency noise, although it simultaneously enhances informative signal components.
In contrast, the corresponding interior data feature $\phi$,
obtained via harmonic extension, enjoys significantly higher regularity: $\phi\in H^{s-2\gamma+\frac{3}{2}}(\Omega)$.
This interior regularity gain compensates for the spectral amplification on the boundary and mitigates the negative effects of noise,
 thereby producing smoother and more structured features that are more favorable for training.
 
\item \textbf{Computation}. Indeed, this pipeline requires the additional computation of $\Phi$ for solving PDEs.
In fact, each learned $\gamma_l$ is associated with its corresponding PDE solver, 
and thus multiple $\{\gamma_l\}$'s demand solving multiple PDEs.
Without a carefully designed differentiable solver, this joint OpL pipeline may increase computational cost.
We shall address this issue in the next section by designing a learnable and differentiable FEM package,
in which all FEM components are implemented using PyTorch’s tensor structures, 
enabling efficient auto-differentiation and batched solutions of multiple PDEs within a unified computational graph.
 \end{itemize}
 
\end{remark}

\section{LA-FEM: the New Learning-Automated FEM Modules in FEALPy}
\label{sec:LAFEM}

To realize the workflow above, we need to integrate FEM and NNs systematically
such that they can push forward and backward each other's output.  

\subsection{Algorithmic Principles of LA-FEM}
Based on \eqref{loss_fun}, we need to compute
\begin{equation}
\label{AD_eq1}
\frac{\partial \mathcal{L}}{\partial \Gamma} = \frac{\partial \mathcal{L}}{\partial \Phi} \cdot \frac{\partial \Phi}{\partial \Gamma}, 
~~~~~ \text{where}~~ K_h(\Gamma) \Phi = F(\Gamma),
\end{equation}
where $F(\Gamma)$ is a matrix of which each column corresponds to one boundary
data pair with the associated $\gamma$,
and then $\Phi$ represents a collection of data feature functions.
For modern DL frameworks such as PyTorch, $\frac{\partial \mathcal{L}}{\partial \Phi}$ can be readily handled by their AD, without users needing to customize the AD behavior through the \texttt{torch.autograd.Function} API.

To our best knowledge, the current PyTorch has NO computing interfaces for $\frac{\partial \Phi}{\partial \Gamma}$, 
i.e., a parametrized FEM solver and its AD with respect to these parameters.
Another concern regarding this workflow is the question of how to solve these $L$ equations in a parallel manner, and this fits naturally to the batched computation tradition in DL.
One key contribution of this work is to develop a novel
Learning-Automated FEM (LA-FEM) module in FEALPy enabling parallel and tensorized computing
using CUDA, AD, and batched computation of PDEs. 
Here we opt to use PyTorch as an example for presentation purpose, yet our implementation works with multiple other backends such as CuPy (numpy) or Jax.

\begin{figure}[h]
    \centering
    \includegraphics[width=0.9\textwidth]{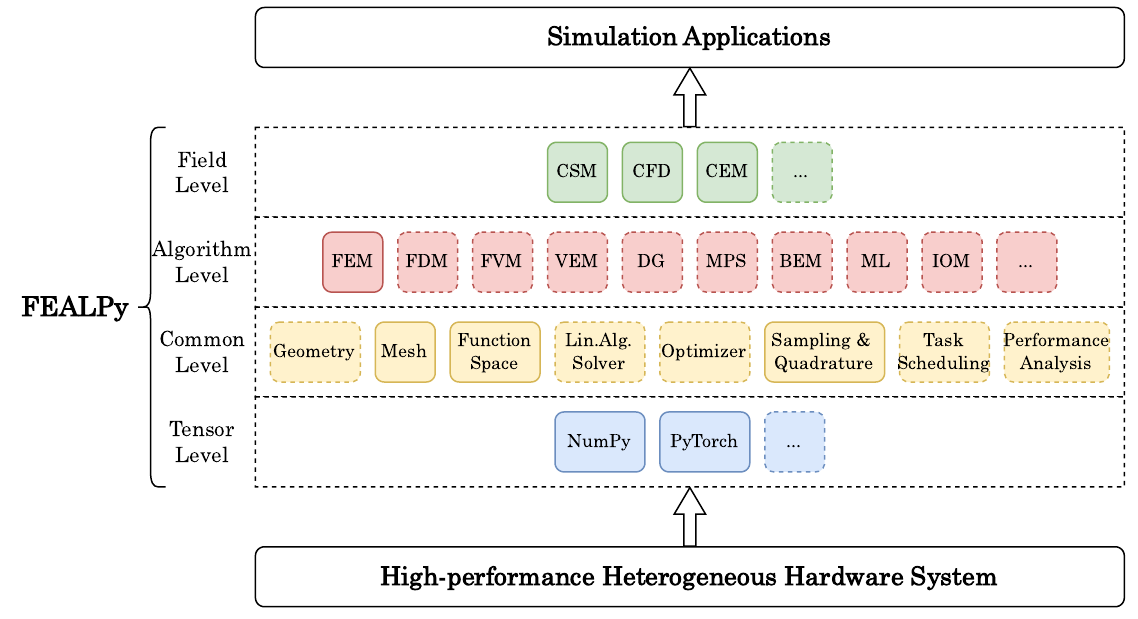}
    \caption{The architecture of FEALPy can be broadly categorized into four levels: 
Tensor level, Common level, Algorithm level, and Field level, 
ranging from low-level to high-level functionality.
}
    \label{fig:fealpy_structure}
\end{figure}

We resort to FEALPy \cite{wei2024fealpy} that is a Python-based open-source library for numerical solutions of PDEs that
incorporates various methods and applications.
The architecture of FEALPy, shown in Fig. \ref{fig:fealpy_structure}, 
already includes an FEM implementation. 
However, previous versions of FEALPy are limited to NumPy, 
which does not have native support for tensor computing with CUDA, AD, and batched computation,
resulting in difficulties in the integration of DL. 
Here, the proposed LA-FEM is designed to incorporate all these features systematically.
Its designing methodology together with some sample usages will be given in Appendix.
The numerical experiments to demonstrate these features are shown in Section \ref{sec:lafem_tests}.

\subsection{Tests for LA-FEM}
\label{sec:lafem_tests}

In this subsection, we shall present a group of numerical experiments to benchmark LA-FEM in tensorized implementation for its scalability in both mesh and batch dimension, as well as efficiency in learning parameters.
For reproducibility, the source codes for the tests of LA-FEM in this section are publicly available in FEALPy~\cite{wei2024fealpy} on GitHub.

\paragraph{Matrix multiplication.} The first test is to verify that the automatic backend manager interface in FEALPY introduces very little overhead.
We refer the reader to Appendix \ref{appendix:backend-manager} for its simple usage and type-agnostic implementation.
We consider a matrix-matrix multiplication of two dense random matrices in the shape of $[N, N]$.
The elements are sampled in float64 type from a uniform distribution in $[0, 1]$.
This matrix multiplication operation (\texttt{matmul}, or \texttt{@}) is executed for 50 times,
and the average wall time is recorded for the NumPy and PyTorch backends in Tab.~\ref{tab:matmul_test}.
The results indicate that our backend manager implementation achieves the expected CUDA acceleration with cuBLAS in PyTorch.

\begin{table}[h]
  \centering
  \begin{tabular}{rrrr}
  \hline
  \multicolumn{1}{c}{\textbf{$N$}} & \textbf{NumPy}(CPU)&\textbf{PyTorch}(CUDA)& \textbf{Speed-up} \\ \hline
  300       & 0.96537     & 0.023918       & 40.4$\times$ \\
  3000      & 172.27      & 0.020323       & 8474.7$\times$ \\
\hline
  \end{tabular}
  \caption{Matrix-matrix multiplication are executed on different backends.
  The results are measured using wall times (in milliseconds), and averaged over 50 runs.}
  \label{tab:matmul_test}
\end{table}

\paragraph{Sparse matrix-vector multiplication.} We conducted a similar experiment on the multiplication of sparse matrices and vectors, which is more commonly encountered in finite element methods and some other approaches.
Sparse matrix-vector multiplication between each square sparse matrix and a random vector was performed,
using \textbf{poisson3Da} (\url{https://sparse.tamu.edu/FEMLAB/poisson3Da}) and
\textbf{cfd1} (\url{https://sparse.tamu.edu/Rothberg/cfd1})
matrices as the test data.
And the elements of the vector are sampled in float64 type from a uniform distribution in $[0,1]$.
The average wall time is recorded for the NumPy and PyTorch backends in Tab.~\ref{tab:sparse_matvec_test}.

\begin{table}[h]
  \centering
\resizebox{\textwidth}{!}
{\begin{tabular}{rrrrrr}
  \hline
  {\textbf{Data}}& \textbf{Size} & \#\textbf{Nonzeros}& \textbf{NumPy}(CPU)&\textbf{PyTorch}(CUDA)& \textbf{Speed-up} \\ \hline
  poisson3Da     & 13,514        & 352,762            & 0.292821           & 0.015521             & 18.9$\times$ \\
  cfd1           & 70,656        & 1,828,364          & 1.320243           & 0.021429             & 61.6$\times$ \\
\hline
  \end{tabular}
}
  \caption{Sparse matrix-vector multiplication are executed on different backends.
  The results are measured using wall times (in milliseconds), and averaged over 50 runs.}
  \label{tab:sparse_matvec_test}
\end{table}

\paragraph{Finite element solvers.} In this test, we benchmark the Numpy and PyTorch backends in LA-FEM using a simple 3D Poisson equation. 
The assembly and solving in PyTorch backend uses GPU, and those in Numpy backend uses CPU. When using PyTorch backend, the extensive usages of \texttt{einsum} in LA-FEM takes the advantage of CUDA parallelism in \texttt{einsum} executions for tensors, and the results show LA-FEM's scalable performance especially for larger tensors with an larger speeding up factor. The results can be seen in Table \ref{torch_exp}. For relatively large problems, 
e.g., more than a million degrees of freedom, the speed-ups of tensorized implementation of FEM solving running on GPUs are drastic over those in NumPy implementation on CPUs, taking fully advantages of running \texttt{einsum} operators on CUDA cores. Note that, the speed-ups of matrix assembly scales with the problem size as well.

\begin{table}[h]
\resizebox{\textwidth}{!}{\begin{tabular}{ccrrrrrrr}
\hline
\multicolumn{1}{c}{\multirow{2}{*}{$h$}} &
  \multirow{2}{*}{\textbf{Order}} &
  \multicolumn{1}{c}{\multirow{2}{*}{\#\textbf{Dofs}}} &
  \multicolumn{2}{c}{\textbf{Wall Time(s) (NumPy)}} &
  \multicolumn{2}{c}{\textbf{Wall Time(s) (PyTorch)}} & 
    \multicolumn{2}{c}{\textbf{Speed-up}}  \\
\cmidrule(lr){4-5}
\cmidrule(lr){6-7}
\cmidrule(lr){8-9}
& & &
  \textbf{Assembly} &
  \textbf{Solve} &
  \textbf{Assembly} &
  \textbf{Solve} &
  \textbf{Assembly} &
  \textbf{Solve}
\\ \hline
$1/128$ & 1 & 16641   & 0.16412 & 0.02102   & 0.75906 & 0.19710  & 0.22$\times$ & 0.11$\times$\\
$1/256$ & 1 & 66049   & 0.45348 & 0.13923   & 0.85730 & 0.30185  & 0.53$\times$ & 0.46$\times$\\
$1/512$ & 1 & 263169  & 1.73856 & 2.20454   & 1.18855 & 0.52189  & 1.46$\times$ & 4.22$\times$\\
$1/512$ & 2 & 1050625 & 4.04250 & 26.56865  & 1.33297 & 3.43641  & 3.03$\times$ & 7.73$\times$\\
$1/512$ & 3 & 2362369 & 8.47439 & 168.72415 & 1.58196 & 15.83130 & 5.36$\times$ & 10.66$\times$\\ \hline
\end{tabular}
}
\caption{Experiments conducted on the same CPU and GPU. A 2D Laplace equation is solved using a uniform triangle mesh on $[0, 1]\times [0, 1]$ with different mesh sizes $h$, and the Lagrange FEM with different orders.
The solvers use Conjugate-Gradient method without preconditioners.
The assembly and solve are measured in seconds using wall times.
FEM on GPU scales nicely as the number of DoFs increases, thanks to the \texttt{einsum}'s automated multi-threaded parallelism.
}
\label{torch_exp}
\end{table}

\begin{figure}[h]
    \centering
    \includegraphics[width=0.6\textwidth]{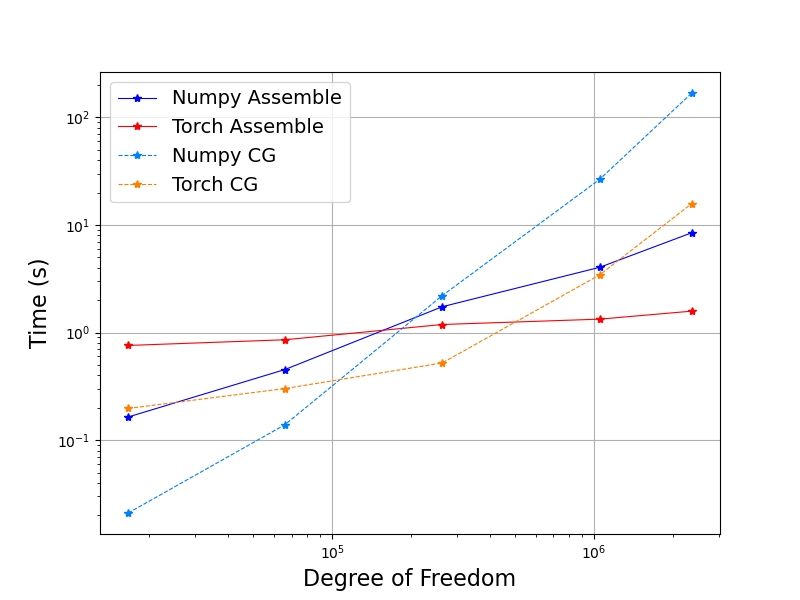}
    \caption{The wall times comparison from Table~\ref{torch_exp}. Both times and the number of Dofs are shown in log-scale.}
\end{figure}

\paragraph{Batched computation.} In this test, the performance of FEM solving using batched computation is compared to that with a serial one on GPUs.
Here we set $L$ Poisson problems on $[0, 1]^3$. The true solutions of them are
\begin{equation}
  u_l = \cos(l\pi x)\cos(l\pi y)\cos(l\pi z), \quad l=1,2,\dots,L,
\end{equation}
and the boundary conditions are set to be of Neumann type with the source term computed correspondingly.
The standard implementation is to assemble and solve the linear system by looping with respect to $L$.  
while in the batched implementation, if memory allocation allows, LA-FEM can solve them in one batch without looping.
Here we set $L=10$, and compared the wall times spent on these two implementations including assembly and solving (average of 10 runs). The results are shown in Table \ref{tab:batch_test}.

\begin{table}[h]
  \centering
  \begin{tabular}{rrrc}
  \hline
  \textbf{DoFs} & \textbf{Serial} & \textbf{Batched} & \textbf{Speed-up} \\ \hline
  35937   & 740.7     & 154.9    & 4.78$\times$ \\
  274625        & 1632.5     & 594.2  & 2.75$\times$ \\ 
\hline
\end{tabular}
  \caption{Wall time (in milliseconds) comparison between the batched computation and looped computation for $L=10$. Both tests are executed on a single NVIDIA A100 80GB GPU. }
  \label{tab:batch_test}
\end{table}

\begin{figure*}[h]
    \begin{center}
    \begin{minipage}{0.48\textwidth}
        \centering
        \includegraphics[width=1\textwidth]{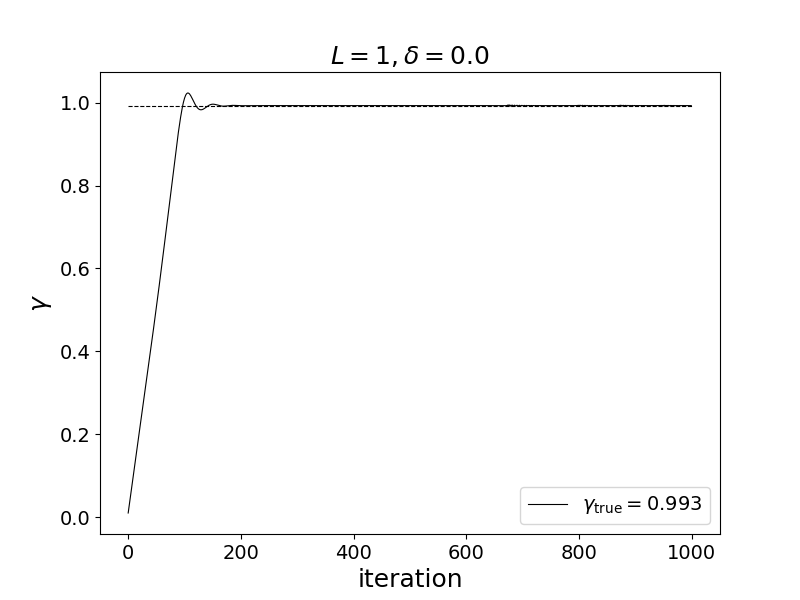}
    \end{minipage}
    ~~~
    \begin{minipage}{0.48\textwidth}
        \centering
        \includegraphics[width=1\textwidth]{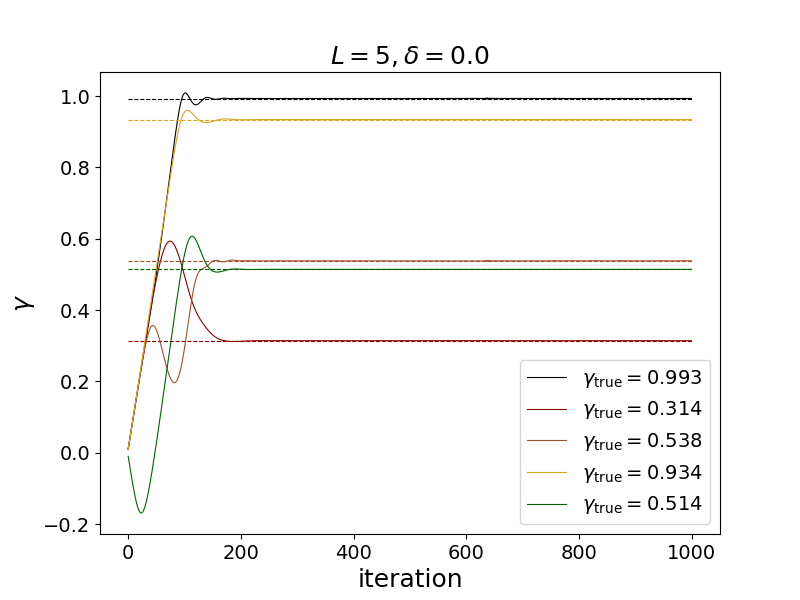}
    \end{minipage}
\end{center}

    \begin{center}
    \begin{minipage}{0.48\textwidth}
        \centering
        \includegraphics[width=1\textwidth]{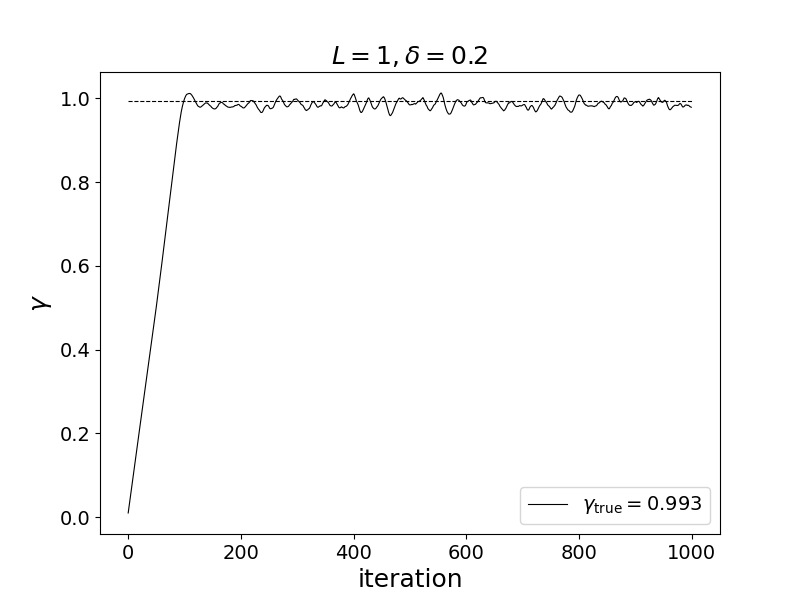}
    \end{minipage}
    ~~~
    \begin{minipage}{0.48\textwidth}
        \centering
        \includegraphics[width=1\textwidth]{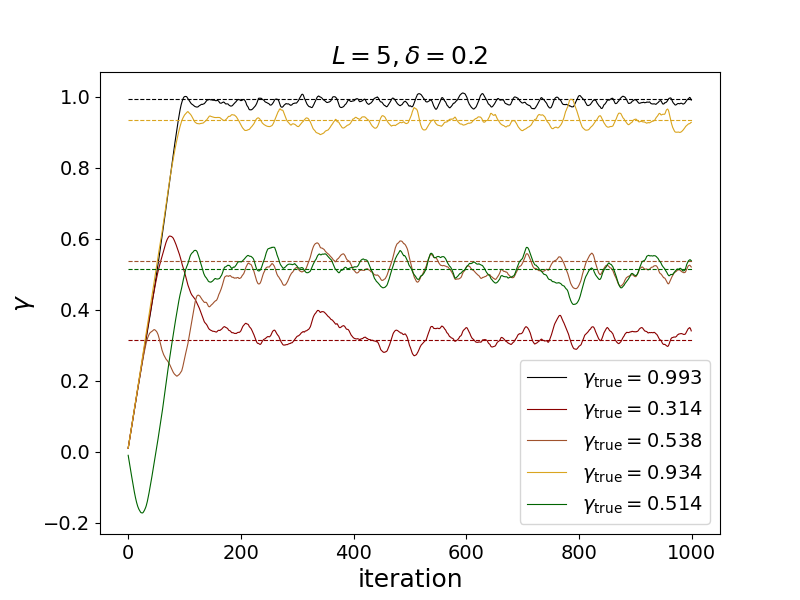}
    \end{minipage}
\end{center}
    \caption{The evolution of $\bgamma$ of the optimization problem \eqref{optim_problem}.
    The native Adam optimizer is used to iterate 1000 steps.
    The dashed line indicates the true value.}
    \label{optim_result}
\end{figure*}

\paragraph{Auto-differentiation in LA-FEM.} In this test, the learning capability of the LA-FEM is benchmarked.
We shall consider a simplified inverse problem to recover multiple fractional orders $\Gamma = (\gamma_1, \gamma_2, \dots, \gamma_L)$ from certain internal data simultaneously.
Let $\Omega_d=[-0.5, 0.5]\times[-0.5, 0.5] \subseteq\Omega$, and for $L$ given $\bgamma^{\textrm{true}}_l$, we consider the solution to \eqref{bc_data_eq2} be $\phi^{\textrm{true}}_l$, $l=1,...,L$, 
corresponding to the boundary data 
$$
\nabla\phi_l\cdot\bfn = \mathfrak{L}^{\gamma_l}_{\partial\Omega}\xi.
$$
Suppose that the true Dirichlet boundary values are accessible on the subdomain $\Omega_{d}$ without evaluating the NtD map,
$$
\phi_D = \sum_{l=1}^L\phi_l^{\text{true}}.
$$
Then, the solving of the following optimization problem is benchmarked:
\begin{align}
\min_{\Gamma = \{ \gamma_1, ..., \gamma_L \} } \quad \mathcal{L}(\Gamma) :=  \left \Vert \sum_{l=1}^L\phi_l - \phi_D \right \Vert^2_{L^2(\Omega_d)}.
\label{optim_problem}
\end{align}
We ran the optimization for 1000 steps on this problem using the native Adam optimizer.
The results in Fig. \ref{optim_result} show that the learned parameter $\Gamma$ converges to $\Gamma^\text{true}$ effectively, for both the cases of
a single and multiple $\gamma$.
We also test the reconstruction when noise appears, and the results are reported in the second row of Fig. \ref{optim_result}.
In this case, even with the presence of a large noise $\delta = 0.2$, approximately $20\%$ error can still be achieved, which is deemed satisfactory within the variance range.
Supported by FEALPy's batched PDE solving capability mentioned in the previous section, the optimization of $\gamma$ is performed in PyTorch's batched tradition without looping through $L$.
In the case of the single $\gamma$ and multiple $\gamma$, the overall optimization took 12 and 23 seconds on 12th Gen Intel® Core™ i7-12700, and took 8 and 10 seconds on NVIDIA A100 80GB GPU, respectively.
See Fig.\ref{tab:learn_gamma_prof} for the detailed profiling results.
Notice that the time required to train the 5-$\gamma$ case is not 5 times as long as the time required to train the single-$\gamma$ case.
\begin{table}[htbp]
\centering
\begin{tabular}{llrrrr}
\hline
\textbf{Device}   & \textbf{L} & \textbf{CPU Time} & \textbf{GPU Time} & \textbf{GFLOPs} & \textbf{Wall Time}
\\ 
\hline \multirow{2}{*}{CPU} & 1 & 12.406    & - & - & 11.897   
\\
                            & 5 & 23.960    & - & - & 23.463
\\ 
\hline
\multirow{2}{*}{GPU}        & 1 & 13.734    & 3.203 & 36.980 & 7.766  
\\
                            & 5 & 14.831    & 8.505 & 17.224 & 9.773    
\\ \hline
\end{tabular}
\caption{The time required to train the 1-$\gamma$ and 5-$\gamma$ cases on CPU and GPU using the PyTorch backend, recorded by \texttt{torch.profiler}. Times are measured in seconds.}
\label{tab:learn_gamma_prof}
\end{table}

\begin{remark}
The proposed method is fundamentally different from PINN 
\cite{2022GaoZahrWang,2022JagtapMaoAdamsKarniadakis,lu2021physics,2019RaissiPerdikarisKarniadakis,2021YangMengKarniadakis,2025WuDuanSunYu}
in the sense that it does not use NNs for approximating PDE solutions.
Instead, the solution is still represented by conventional FE functions; only some unknown parameters are included as unknowns for training through backpropagation.
In the present work, only fractional orders are used as unknowns.
Certainly, for more complicated problems, NNs can be also coupled with FEM in the LA-FEM framework,
which will be an interesting direction in the future.
\end{remark}

\section{Numerical experiments}
\label{sec:num}

In this section, the proposed $\gamma$-deepDSM is benchmarked to test its reconstruction capability for multiple inclusions, in which the role of  learning $\gamma$ is demonstrated.
Source codes of FEALPy are publicly available at \url{https://github.com/weihuayi/fealpy}. The reproducible tests based on LA-FEM are available at
\url{https://github.com/weihuayi/fealpy/tree/master/app/lafem-eit} as an app based on FEALPy.

The computational domain is set to be $\Omega=(-1, 1)\times(-1, 1)$.
The incident current data $\{g_{N,l}\}$ on the boundary are chosen as
$$
g_{N,l}(\bfx)=\cos(l\theta(\bfx)), \quad \bfx\in\partial\Omega, \quad l=1,2,3,4,5,6,8,16,
$$
and when being discretized, the instances are concatenated as the channel dimension for the NN.
$\theta(\bfx)$ denotes the polar angle of $\bfx$ on $\partial\Omega$.
For each inclusion $\sigma$, solving \eqref{eq_eit} with $g_{N,l}$ FEM on triangular meshes shall yield an approximation to $g_{D,l}$ by restricting the numerical solution on the boundary.
When obtaining $g_{D,l}$, noises are added to the FEM approximation to mimic the measurement error in practice.
Let $G(x)$ follow a Gaussian distribution, 
and let $G(x)$ and $G(y)$ be independent if $x\neq y$, i.e., the Gaussian noise.
We consider a noise in the form of 
$u^{\delta} = u + \delta(-\Delta_{\partial \Omega})^{-0.75}(Gu)$
with $\delta$ being the parameter to control the intensity, which is referred to as the noise level in subsequent tests.
The noisy data in both present in training set and test set.
Since in the real scenario the incident current $g_N$ is generated by the instrument, we assume that measurement noise for $g_N$ is negligible.
The following three settings are tested and compared:

\begin{itemize} 
\item [(\textbf{S1})] \label{exp:S1} $\gamma=0$. The fLB order is fixed to zero, i.e., no LB operator is applied. 
\item [(\textbf{S2})] \label{exp:S2}There is a single learnable $\gamma$ as the order of fLB operator applied on the scattered data, and initialized as zero.
\item [(\textbf{S3})] \label{exp:S3} There are multiple learnable $\gamma$'s, and initialized as zeros.
\end{itemize}

We have also tested the random initial guess for $\gamma$ but observe no difference from the zero initial guess.
Thus, the following experiments are all based on the zero initial.
In each case, the injected noise has three levels: $\delta$=0\%, 1\% and 5\%.
The pointwise binary cross entropy is chosen as the Loss function.
For fairness and simplicity, a vanilla SGD is chosen as the optimizer with a learning rate $10^{-3}$, momentum $0.9$, and weight decay $10^{-9}$.
The training set consists of $10000$ samples, while the test set consists of $2000$ samples.

\subsection{Basic results}

We first consider the case that $\Omega$ contains two medium with different conductivity $10$ (inclusion) and $1$ (background).
The distribution of random inclusion is generated by the following procedure:
\begin{enumerate}
  \item Sampling $N_c$ 2-D circle centers $(c_{ix}, c_{iy})$ from uniform distribution $U(-0.8, 0.8), i=1,2,\dots,N_c$.
  \item Sampling radius $r_i\sim U(0.1, b_i), b_i=\min\{0.9-|c_{ix}|, 0.9-|c_{iy}|\}, i=1,2,\dots,N_c$.
  \item The inclusion $D$ is defined as the union of these $N_c$ circles.
\end{enumerate}
We refer readers to Table \ref{validate_sets} for the notation of various test data sets, 
where \textbf{CIR3} is an in-distribution dataset for training models in this subsection.

\begin{figure*}[htbp]
    \begin{center}
    \begin{minipage}{0.48\textwidth}
        \centering
        \includegraphics[width=1.0\textwidth]{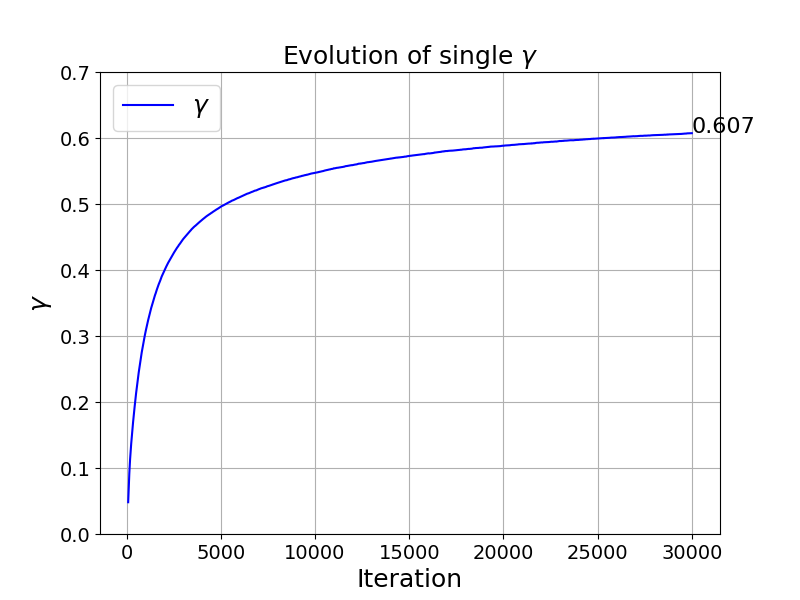}
    \end{minipage}
    \begin{minipage}{0.48\textwidth}
        \centering
        \includegraphics[width=1.0\textwidth]{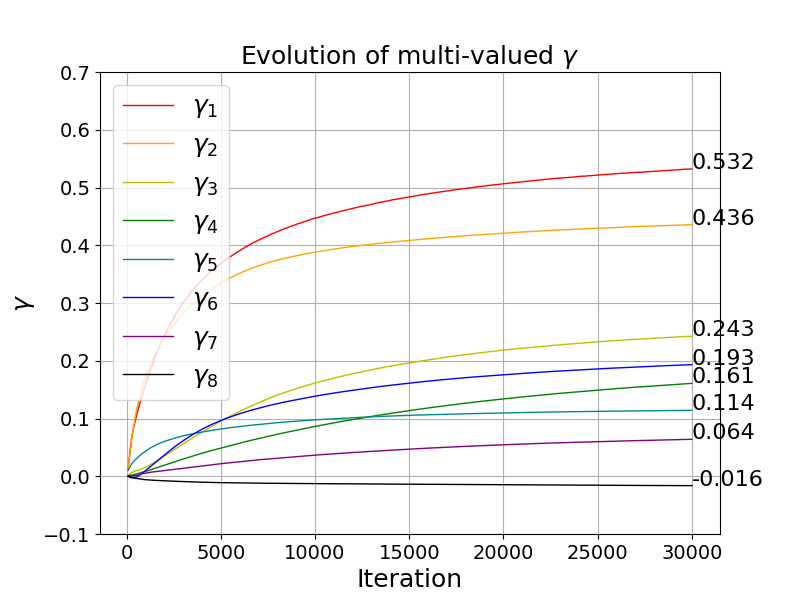}
    \end{minipage}
    \end{center}
    \caption{The training dynamics of $\gamma$ for the single $\gamma$ and multiple $\gamma$'s. Here the fLB order $\gamma$ is optimized like neural network parameters, and $\{\gamma_i\}_{i=1}^8$ represent the values corresponding to $g_{N,l_i}$, $l_i=1,2,3,4,5,6,8,16$.
    It is observed that the higher incident frequency results a lower value of the learned $\gamma$.}
    \label{evo_s}
\end{figure*}

Now we present training and test results for the setting described in \hyperref[exp:S1]{(\textbf{S1})}-\hyperref[exp:S1]{(\textbf{S3})}.
For both the single-$\gamma$ case and multi-$\gamma$ case, it is observed that these $\gamma$'s converge to different values, which is likely attributed to the variation in the boundary data and the noise. 
The results are reported in Fig. \ref{evo_s}.
We observe an inverse correlation between the learned fractional order and the incident frequency; namely, lower-frequency boundary inputs tend to produce larger learned values of $\gamma$. This behavior is consistent with the role of $\gamma$ as a spectral filter: while it amplifies informative high-frequency signal components, it also enhances noise, and therefore must strike a balance to optimize reconstruction accuracy. Although a fully rigorous theoretical characterization remains challenging, we provide an elementary analytical explanation of this phenomenon in Section~\ref{subsec:fLBanalysis}.

\begin{figure*}[htbp]
    \begin{center}
    \begin{minipage}{0.48\textwidth}
        \centering
        \includegraphics[width=1.0\textwidth]{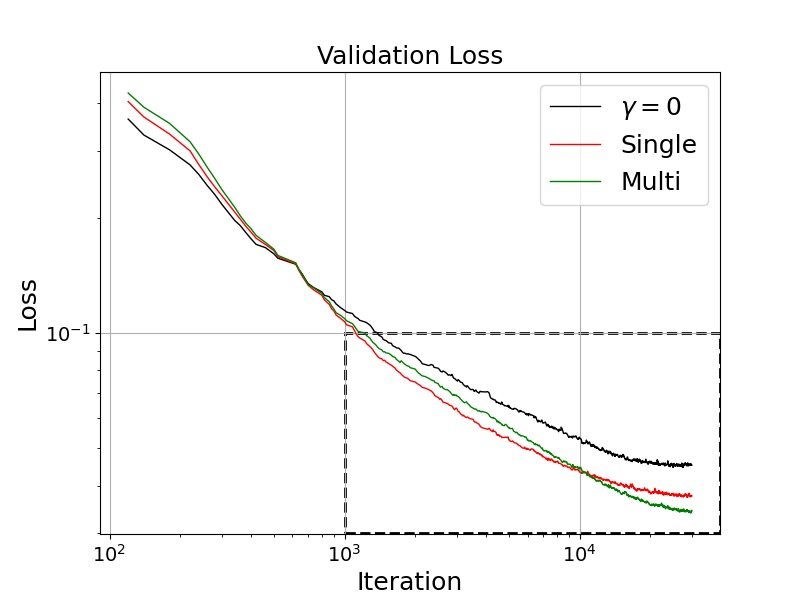}
    \end{minipage}
    \begin{minipage}{0.48\textwidth}
        \centering
        \includegraphics[width=1.0\textwidth]{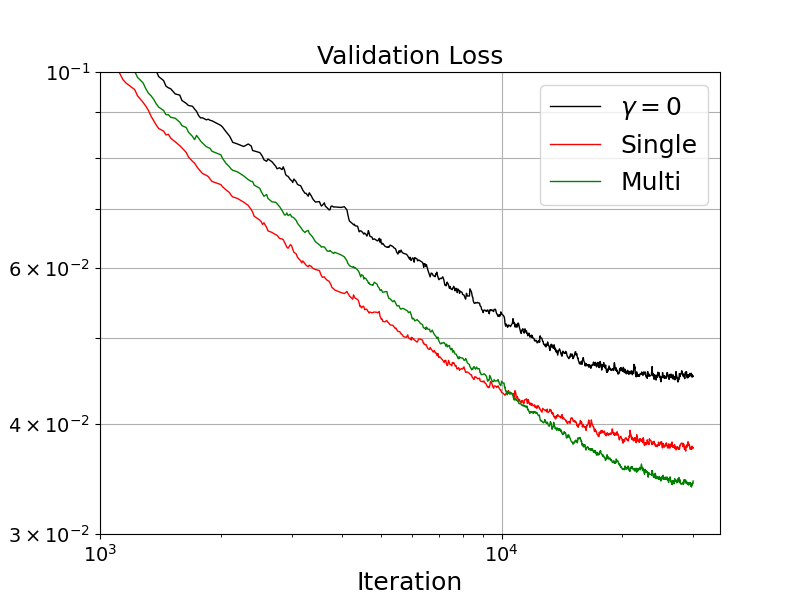}
    \end{minipage}
    \end{center}
    \caption{Convergence of the validation losses. 
    The left plot shows the dynamics over all iterations, 
    and the right one zooms in near $10^4$ iterations.  
    Both the single-$\gamma$ and multi-$\gamma$ models show at least 15\% improvement over the non-learnable $\gamma=0$ one.}
    \label{evo_nn}
\end{figure*}

\begin{figure*}[htbp]
  \begin{center}
    \begin{minipage}{0.48\textwidth}
        \centering
        \includegraphics[width=1.0\textwidth]{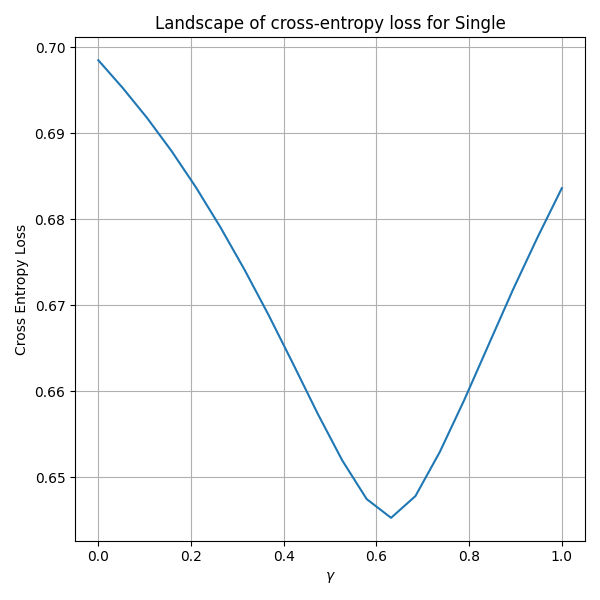}
    \end{minipage}
    \begin{minipage}{0.48\textwidth}
        \centering
        \includegraphics[width=1.0\textwidth]{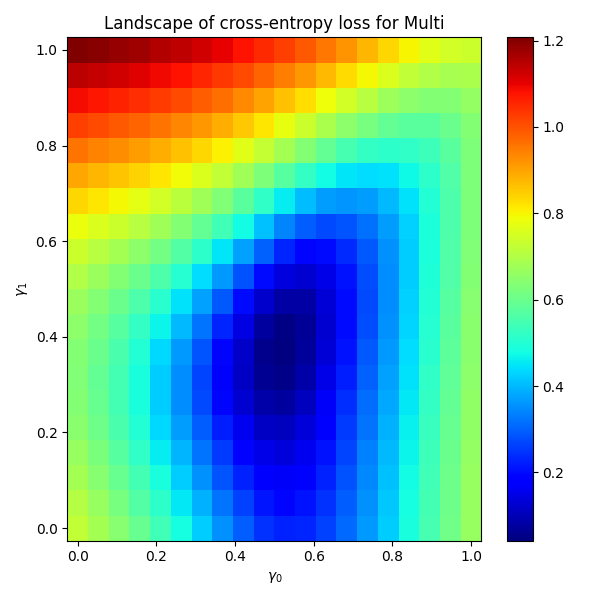}
    \end{minipage}
  \end{center}
  \caption{The landscape of the loss function for the single $\gamma$ and multiple $\gamma$'s.}
  \label{fig:landscape}
\end{figure*}

The plot of loss function values versus training iterations is shown in Fig. \ref{evo_nn} for the noiseless case ($\delta =0\%$).
The multi-$\gamma$ model does not demonstrate a smaller validation loss than the single-$\gamma$ one until $10^4$ steps, 
but ultimately surpassed the single model as the training continues.
It may be attributed to that multiple $\gamma$'s make the landscape more complex and thus increases the optimization complexity. 
We refer readers to Fig. \ref{fig:landscape} for the landscape of the loss function with respect to $\gamma$. 
We employ the Adam optimizer during training, which treats the parameters in LA-FEM and those in the NN uniformly. 
While convenient, this unified treatment may not be the most efficient strategy for such a coupled FEM-NN architecture. 
Developing more advanced or adaptive optimization schemes for such architecture would be an interesting direction for future research.

\begin{figure*}[htbp]
  \begin{center}
  \begin{minipage}{0.48\textwidth}
    \centering
    \includegraphics[width=1.0\textwidth]{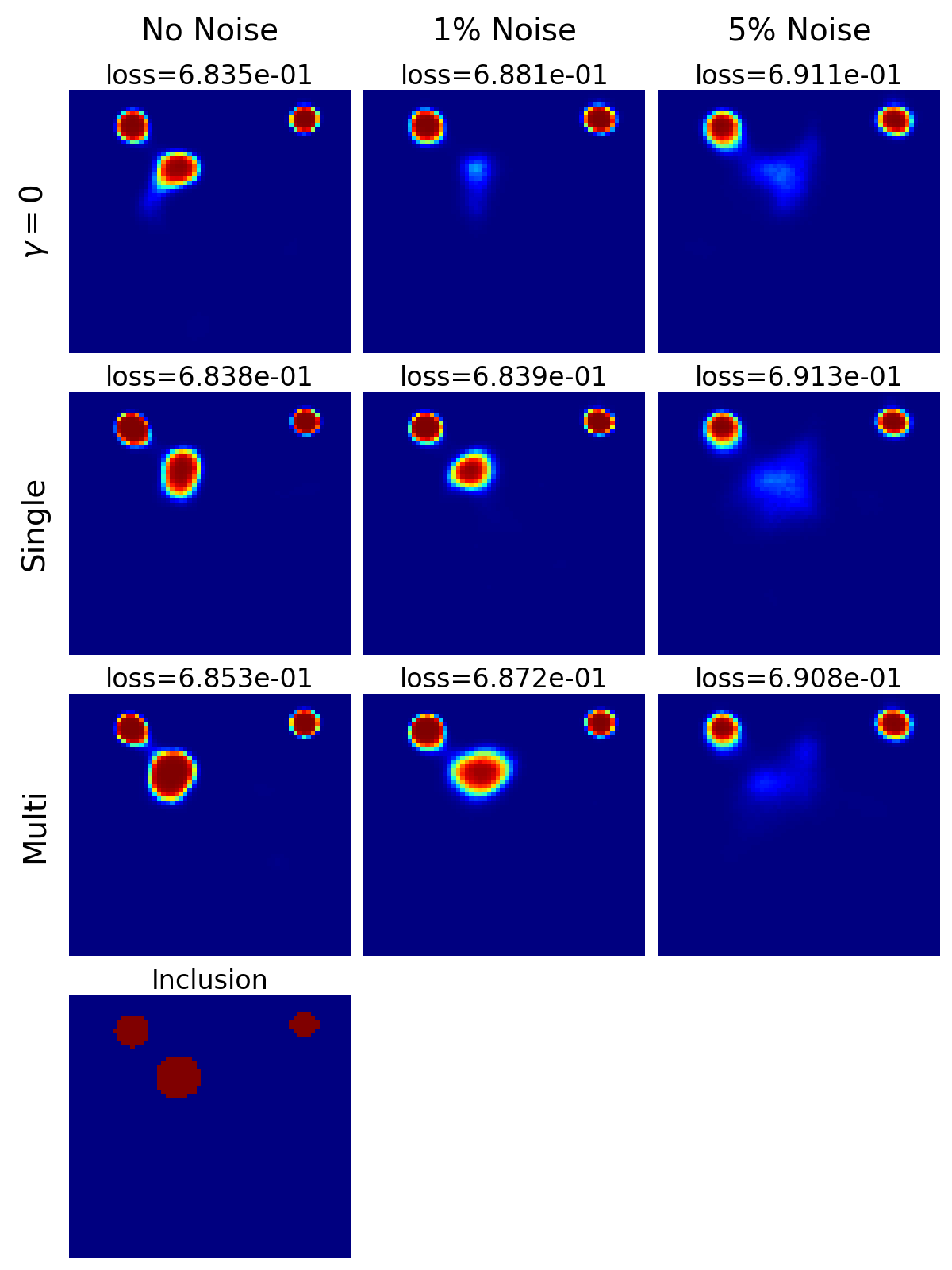}
  \end{minipage}
  \begin{minipage}{0.48\textwidth}
    \centering
    \includegraphics[width=1.0\textwidth]{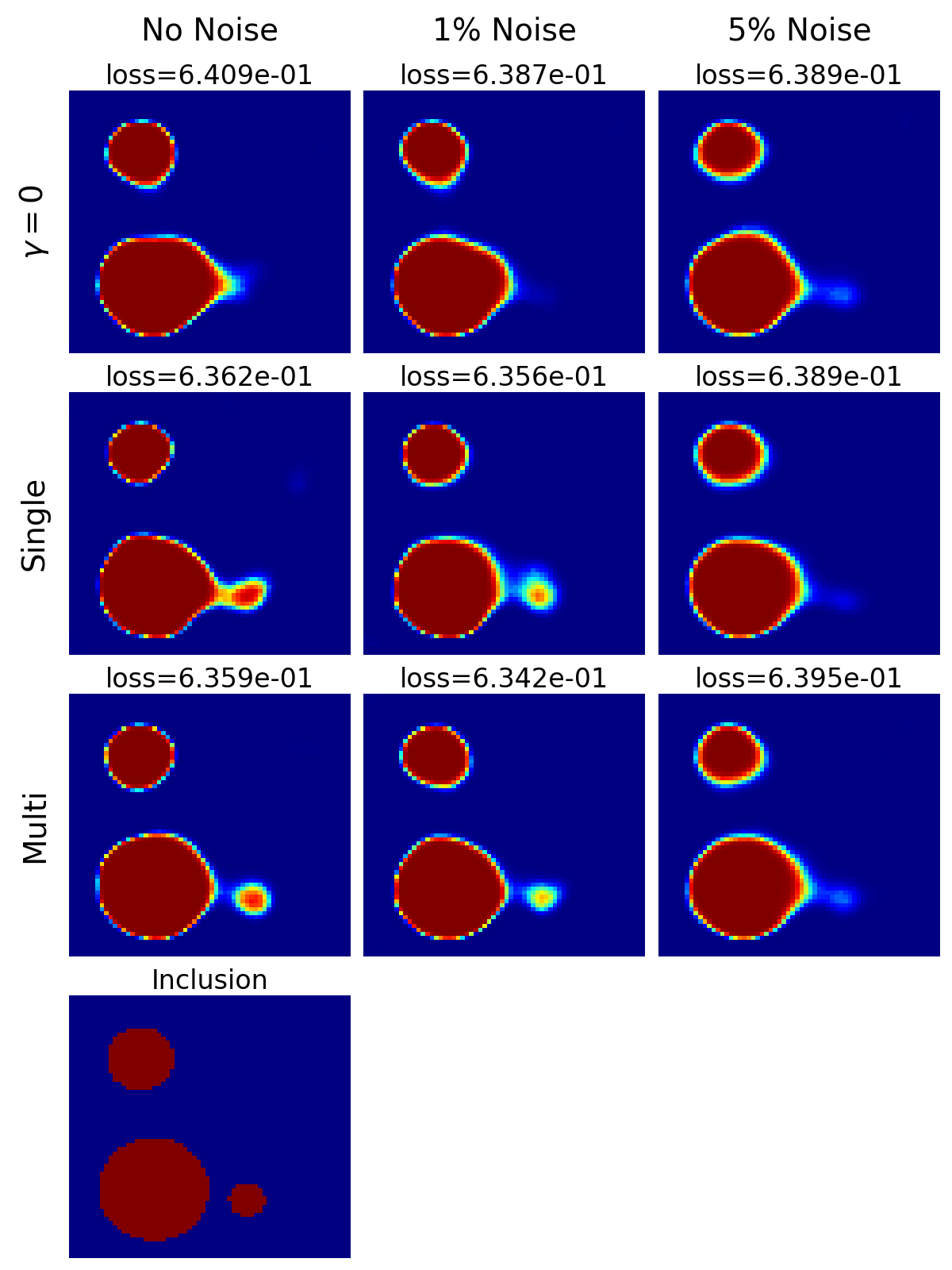}
  \end{minipage}
  \begin{minipage}{0.48\textwidth}
    \centering
    \includegraphics[width=1.0\textwidth]{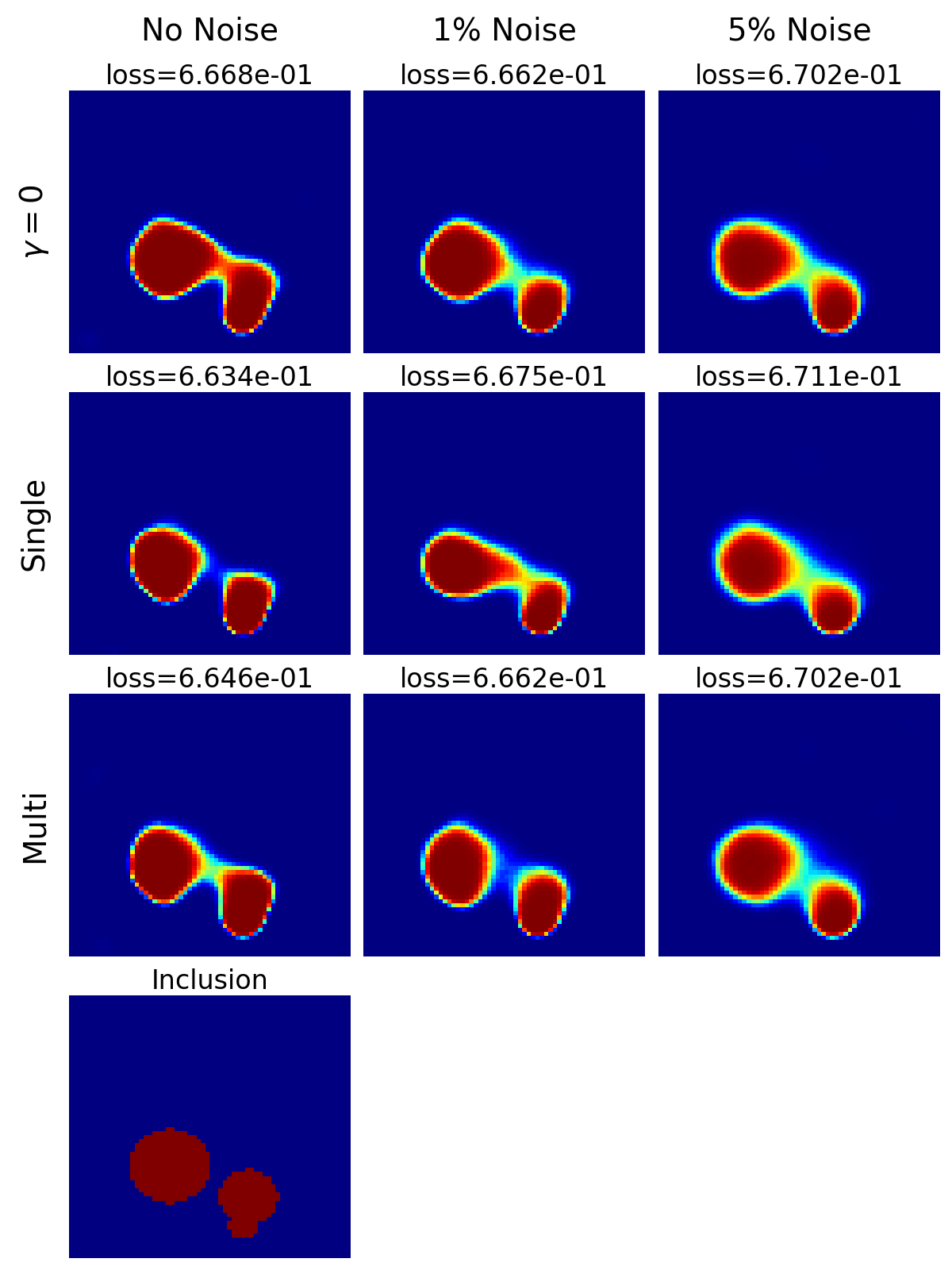}
  \end{minipage}
  \begin{minipage}{0.48\textwidth}
    \centering
    \includegraphics[width=1.0\textwidth]{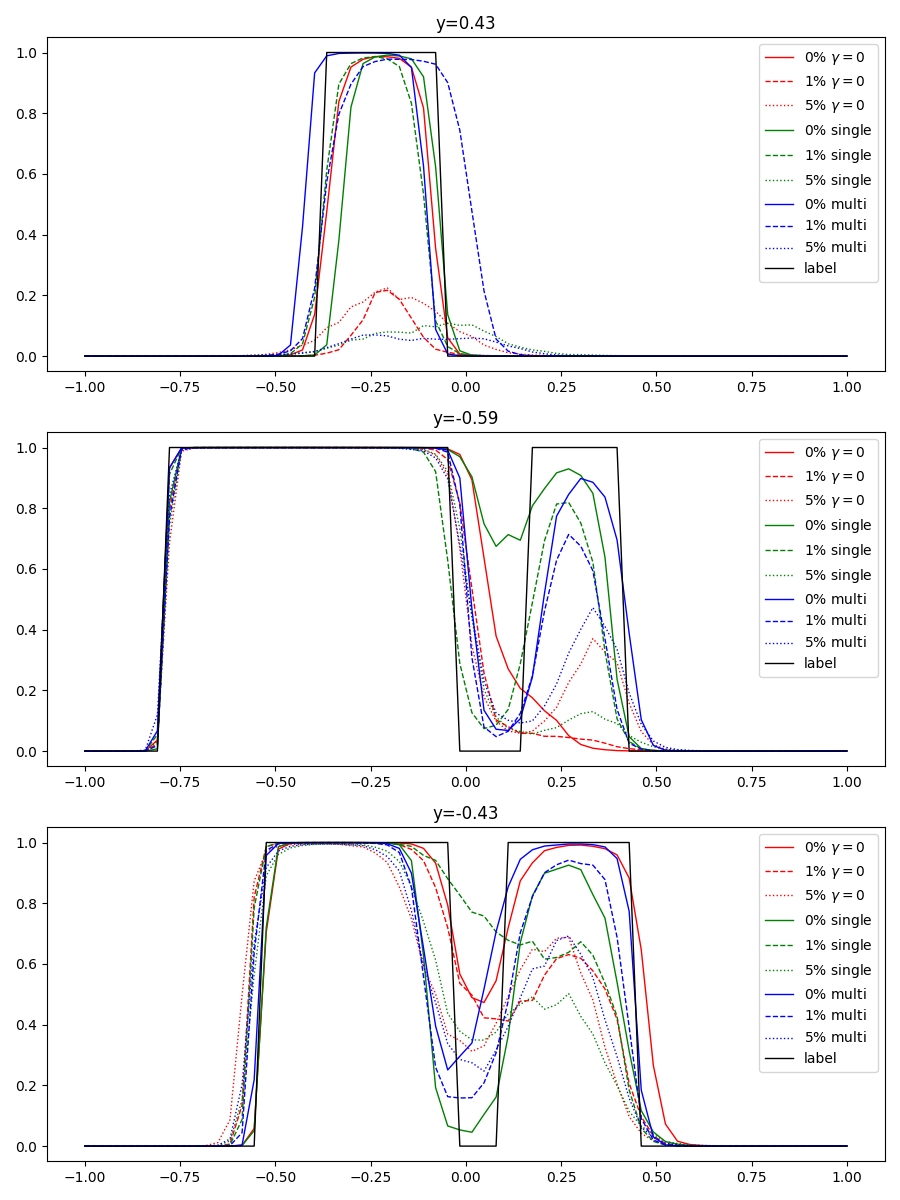}
  \end{minipage}
  \end{center}
  \caption{The in-distribution reconstructed images produced by the models of
  ``$\gamma=0$'' (top), a single learnable $\gamma$ (middle) and multiple
  learnable $\gamma$'s (bottom).
  The white dotted circles denote the ground truth.
  The last axis shows the cross-section of the reconstructed images at a specific y coordinate in order.}
  \label{vis_cir3}
\end{figure*}

All models are trained for 5 times with different random seeds, and the mean and
relative standard deviation (RSD) of the test loss are reported in Table~\ref{out_of_dist}
and Table~\ref{out_of_dist_dice}, for the binary cross entropy (BCE) and dice loss, respectively.
Here, the column associated with \textbf{CIR3} presents an in-distribution data set.
In the noiseless case, the validation loss in the trainable single-$\gamma$ model is about 23\% lower than $\gamma=0$, 
and the multi-$\gamma$ model is about 25\% less than the single-$\gamma$.
In particular, we refer readers to Fig. \ref{vis_cir3} for the case that the small circle and the center circle is likely to be ignored or missed by the non-learnable-$\gamma$ model ($\gamma=0$), 
especially with the presence of noise. 
Remarkably, $\gamma$-deepDSM can visibly better locate these small inclusions and yield sharper boundaries.

\begin{figure*}[htbp]
    \begin{center}
    \begin{minipage}{0.55\textwidth}
        \includegraphics[width=1.0\textwidth]{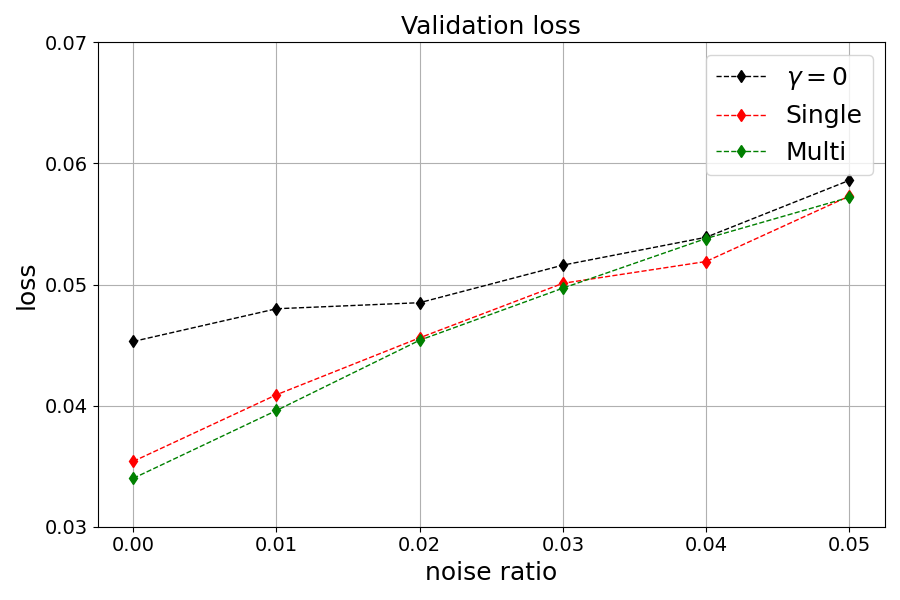}
    \end{minipage}
    \end{center}
    \caption{Validation loss versus the increasing noise level.}
    \label{ter}
\end{figure*}

\begin{figure*}[htbp]
  \begin{center}
    \begin{minipage}{0.48\textwidth}
      \centering
      \includegraphics[width=1.0\textwidth]{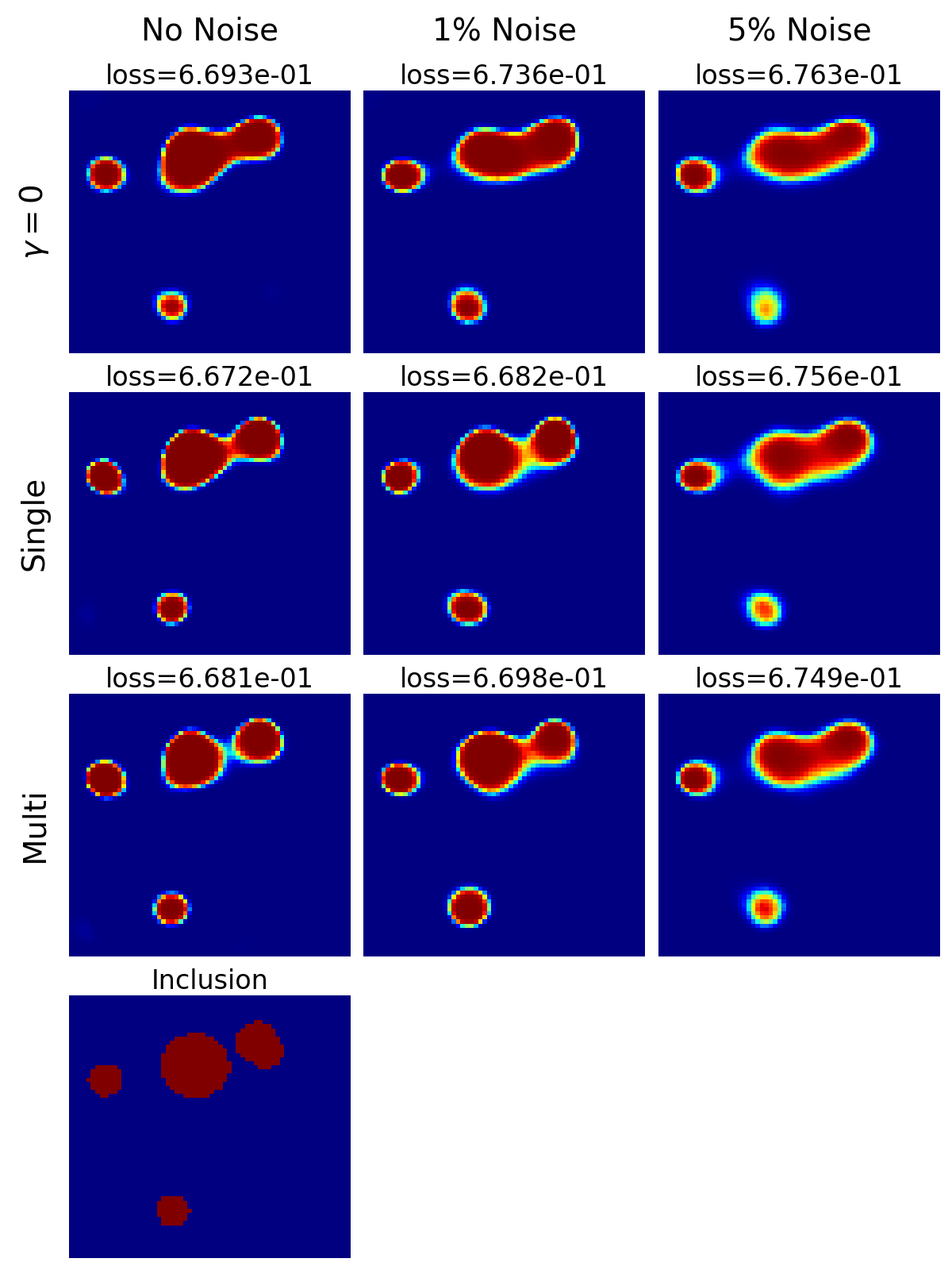}
    \end{minipage}
    \begin{minipage}{0.48\textwidth}
      \centering
      \includegraphics[width=1.0\textwidth]{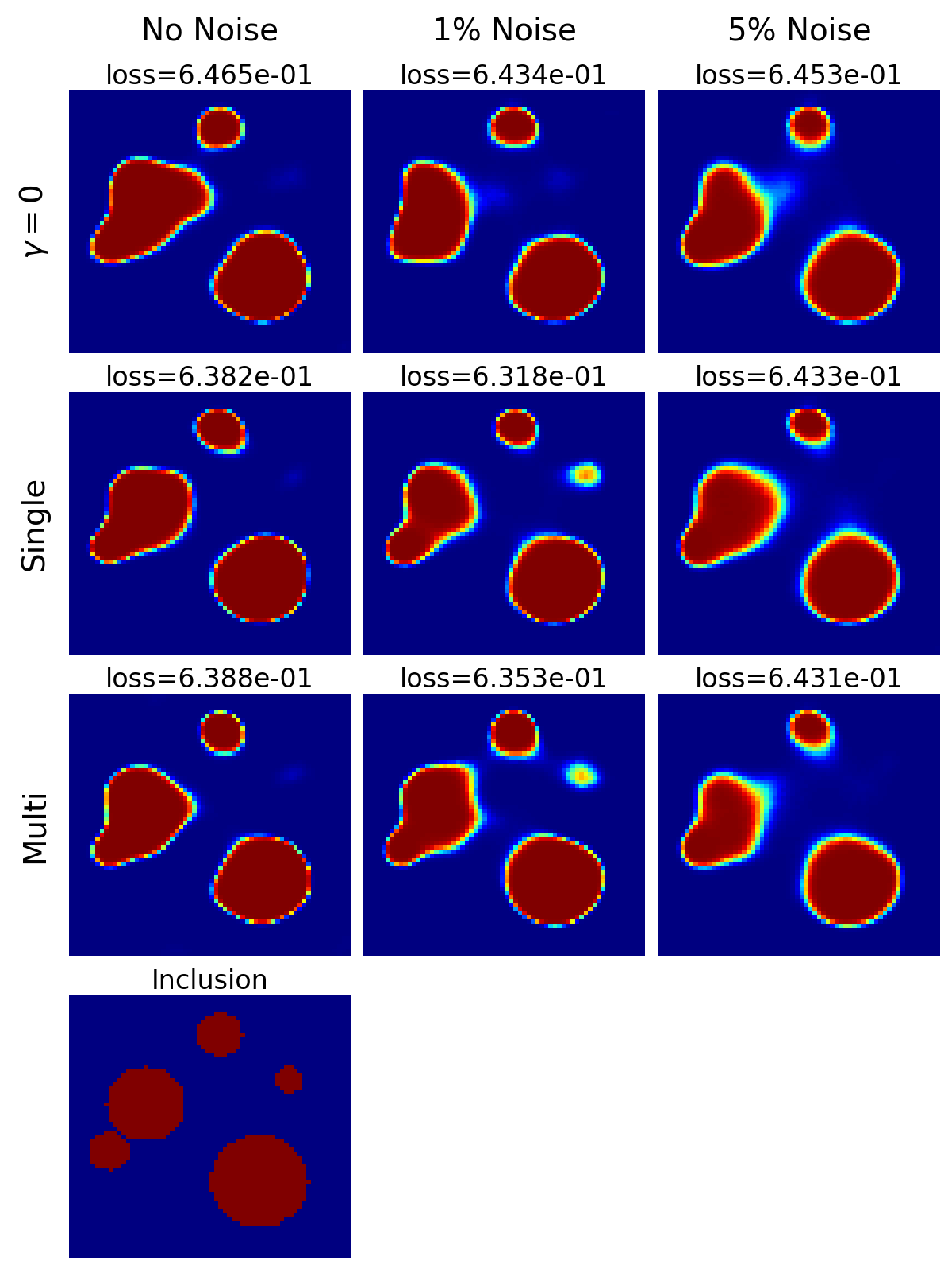}
    \end{minipage}
    \begin{minipage}{0.48\textwidth}
      \centering
      \includegraphics[width=1.0\textwidth]{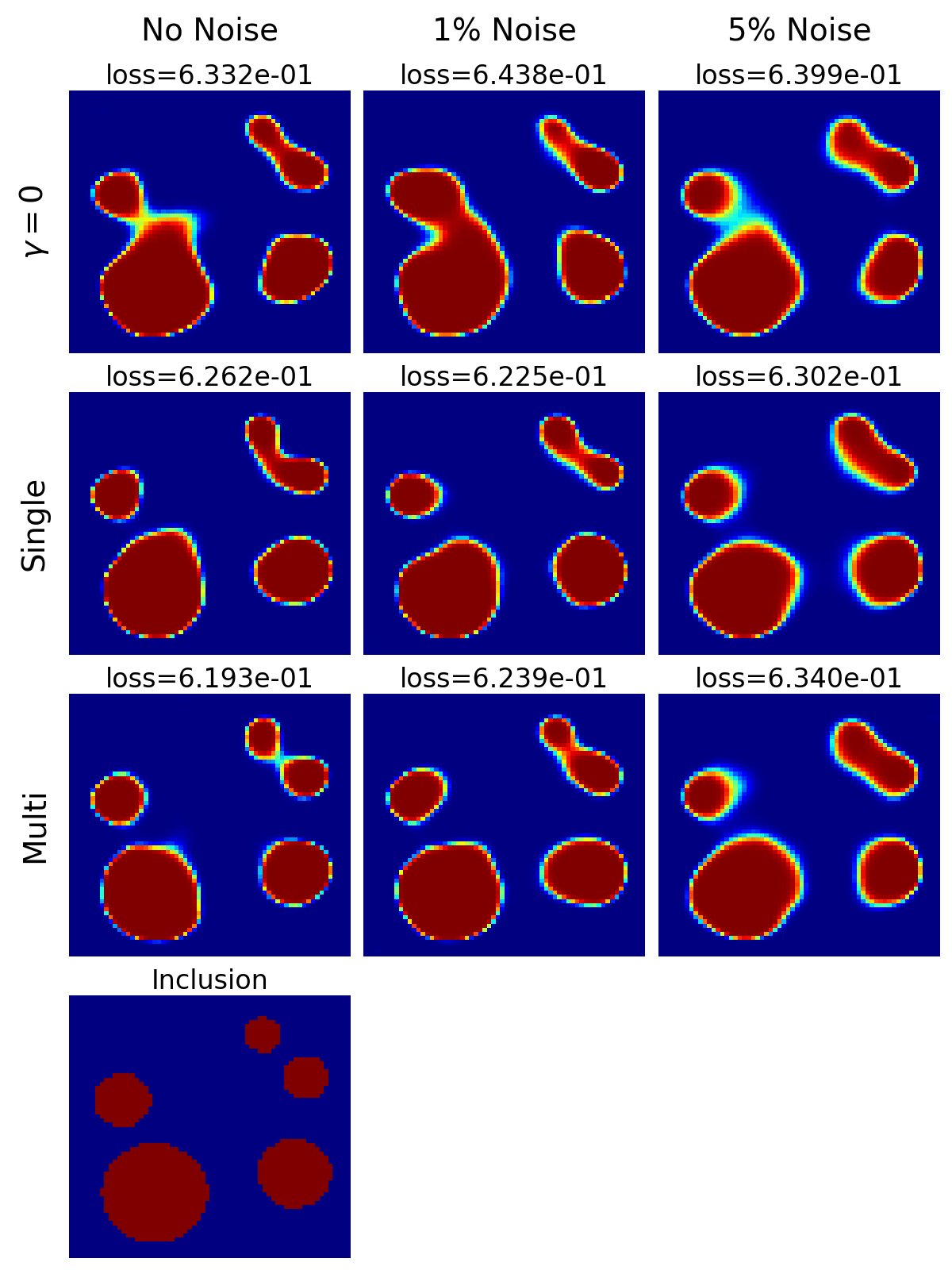}
    \end{minipage}
    \begin{minipage}{0.48\textwidth}
      \centering
      \includegraphics[width=1.0\textwidth]{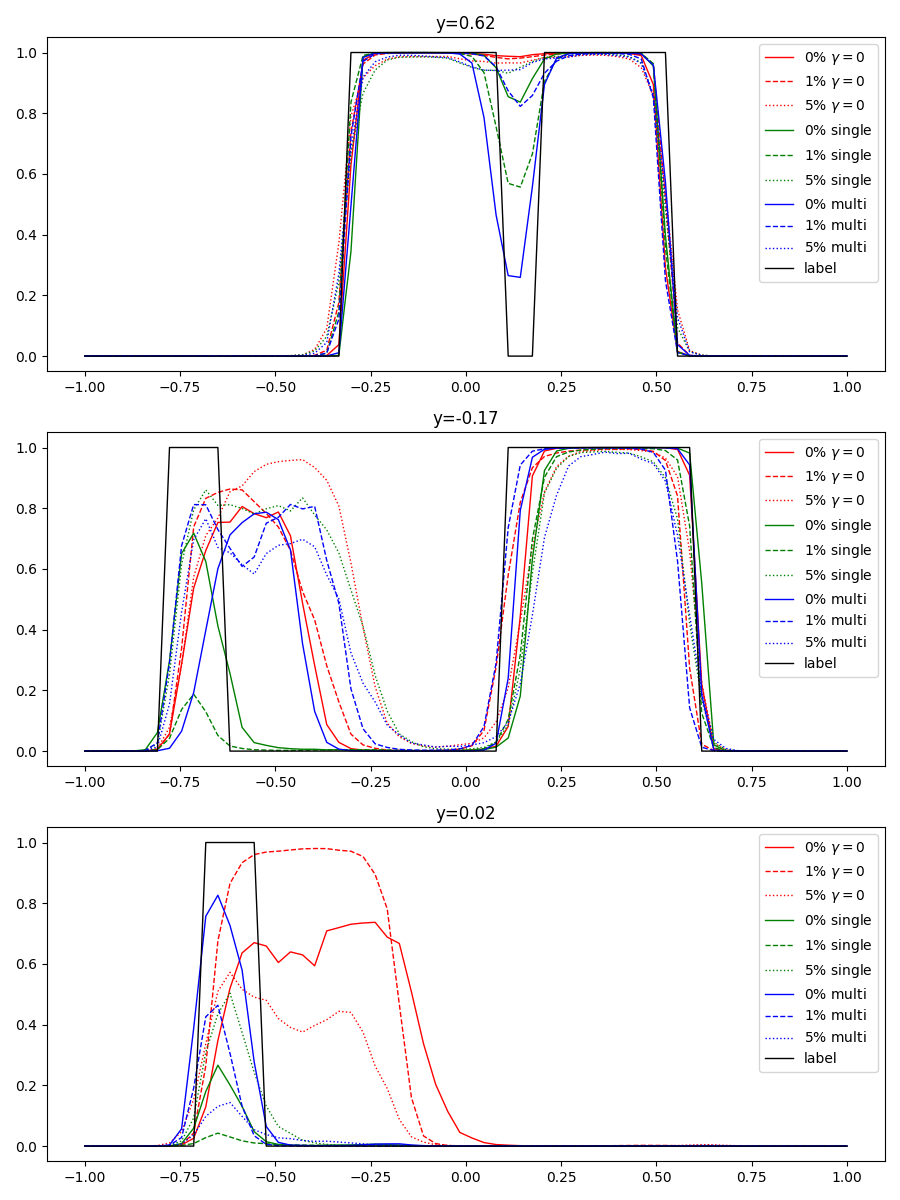}
    \end{minipage}
  \end{center}
  \caption{The out-of-distribution reconstructed images produced by the models of ``$\gamma=0$'' (top), a single learnable $\gamma$ (middle) and multiple learnable $\gamma$'s (bottom).
  The white dotted circles denote the ground truth.
  The last axis shows the cross-section of the reconstructed images at a specific y coordinate in order.}
  \label{vis_cir5}
\end{figure*}

Next, we perform out-of-distribution (OOD) tests where the validation data set has a significantly different distribution than the training set.
The inclusions in the training dataset are fixed with $N_c=3$ circles.
Then, the different number of circles ($N_c=2, 4, 5$) are used to test the generalization ability.
Also, all models are trained for 5 times with different random seeds, and the mean and
relative standard deviation (RSD) of the test loss are reported in Table~\ref{out_of_dist}
and Table~\ref{out_of_dist_dice}, for the BCE and dice loss, respectively.
From these results, we can also see that the learnable fLB orders can achieve almost $\mathbf{15\%}$ improvement compared with the model with non-learnable-$\gamma$.
For an OOD example, we refer to readers for five circular inclusions in Fig. \ref{vis_cir5}.
Therein, even though the non-learnable-$\gamma$ ($\gamma=0$) can locate the approximate inclusion locations, 
it fails to visibly separate inclusion boundaries, especially when two inclusions are relatively close to each other; see the first and the third cases in Fig.~\ref{vis_cir5}.
In the meantime, it may ignore small inclusions; see the second case in Fig.~\ref{vis_cir5} for instance. 
All the reconstructions in these scenarios are largely improved by making $\gamma$ learnable in $\gamma$-deepDSM.
\begin{table}[htbp]
\centering
\caption{Notation and setup of different datasets.}
\label{validate_sets}
\begin{tabular}{cccccc}
Notation & Size                 & Channels           & Mesh Size                      & Sigma                   & Number of Circles \\ \hline
CIR2                    & \multirow{4}{*}{2000} & \multirow{4}{*}{8} & \multirow{4}{*}{$64\times 64$} & \multirow{4}{*}{$10,1$} & 2                 \\
CIR3 &  &  &  &  & 3 \\
CIR4 &  &  &  &  & 4 \\
CIR5 &  &  &  &  & 5 \\ \hline
\end{tabular}
\end{table}

\begin{table}[htbp]
\centering
\caption{Test loss (binary cross entropy) on the OOD datasets (CIR3--CIR5),
compared to the in-distribution data set.
The reported mean and RSD are calculated based on 5 runs with different random seeds.}
\label{out_of_dist}
\small
\begin{tabular}{ll|ll|ll|ll|ll}
\hline
$\delta$       & \textbf{Model}  & \textbf{CIR2}    &          & \textbf{CIR3}    &          & \textbf{CIR4}    &          & \textbf{CIR5}    &          \\
               &                 & mean             & RSD      & mean             & RSD      & mean             & RSD      & mean             & RSD      \\ \hline
0\%            & $\gamma=0$      & 0.02580          & 1.020\%  & 0.04530          & 1.867\%  & 0.07434          & 3.171\%  & 0.10279          & 3.504\%  \\
               & Single-$\gamma$ & 0.01945          & 1.288\%  & 0.03553          & 1.926\%  & 0.05744          & 3.269\%  & 0.07972          & 3.763\%  \\
               & \textit{Improve} & 24.61\%          &          & 21.57\%          &          & 22.73\%          &          & 22.44\%          &          \\
               & Multi-$\gamma$  & \textbf{0.01903} & 1.796\%  & \textbf{0.03456} & 2.247\%  & \textbf{0.05543} & 2.894\%  & \textbf{0.07710} & 3.041\%  \\
               & \textit{Improve} & 26.24\%          &          & 23.71\%          &          & 25.44\%          &          & 24.99\%          &          \\ \hline
1\%            & $\gamma=0$      & 0.02942          & 1.184\%  & 0.04826          & 1.555\%  & 0.07590          & 1.917\%  & 0.10272          & 2.165\%  \\
               & Single-$\gamma$ & 0.02529          & 1.026\%  & 0.04112          & 0.392\%  & 0.06142          & 1.288\%  & 0.08192          & 1.559\%  \\
               & \textit{Improve} & 11.90\%          &          & 12.54\%          &          & 19.08\%          &          & 20.52\%          &          \\
               & Multi-$\gamma$  & \textbf{0.02444} & 1.037\%  & \textbf{0.04012} & 0.801\%  & \textbf{0.06031} & 1.505\%  & \textbf{0.08079} & 2.233\%  \\
               & \textit{Improve} & 16.93\%          &          & 16.87\%          &          & 20.54\%          &          & 21.35\%          &          \\ \hline
5\%            & $\gamma=0$      & 0.03892          & 1.557\%  & 0.05895          & 1.499\%  & 0.08471          & 1.126\%  & 0.10881          & 1.386\%  \\
               & Single-$\gamma$ & 0.03918          & 0.839\%  & 0.05818          & 1.496\%  & 0.08234          & 2.233\%  & 0.10411          & 3.165\%  \\
               & \textit{Improve} & -0.67\%          &          & 1.31\%           &          & 2.80\%           &          & 4.32\%           &          \\
               & Multi-$\gamma$  & \textbf{0.03873} & 2.378\%  & \textbf{0.05775} & 1.624\%  & \textbf{0.08181} & 1.846\%  & \textbf{0.10345} & 2.246\%  \\
               & \textit{Improve} & 0.49\%           &          & 2.04\%           &          & 3.42\%           &          & 4.93\%           &          \\ \hline
\end{tabular}
\end{table}

\begin{table}[htbp]
\centering
\caption{Test loss (dice) on the OOD datasets (CIR3--CIR5),
compared to the in-distribution data set.
The reported mean and RSD are calculated based on 5 runs with different random seeds.}
\label{out_of_dist_dice}
\small
\begin{tabular}{ll|ll|ll|ll|ll}
\hline
$\delta$       & \textbf{Model}  & \textbf{CIR2}    &          & \textbf{CIR3}    &          & \textbf{CIR4}    &          & \textbf{CIR5}    &          \\
               &                 & mean             & RSD      & mean             & RSD      & mean             & RSD      & mean             & RSD      \\ \hline
0\%            & $\gamma=0$      & 0.11965          & 0.873\%  & 0.11131          & 1.091\%  & 0.11057          & 1.136\%  & 0.11297          & 1.073\%  \\
               & Single-$\gamma$ & 0.08890          & 1.252\%  & 0.08721          & 1.298\%  & 0.08821          & 1.310\%  & 0.09030          & 1.109\%  \\
               & \textit{Improve} & 25.70\%          &          & 21.65\%          &          & 20.22\%          &          & 20.07\%          &          \\
               & Multi-$\gamma$  & \textbf{0.08847} & 2.011\%  & \textbf{0.08678} & 1.699\%  & \textbf{0.08726} & 1.333\%  & \textbf{0.08958} & 1.424\%  \\
               & \textit{Improve} & 26.06\%          &          & 22.04\%          &          & 21.08\%          &          & 20.70\%          &          \\ \hline
1\%            & $\gamma=0$      & 0.13888          & 1.036\%  & 0.12416          & 0.953\%  & 0.11965          & 0.984\%  & 0.11993          & 1.006\%  \\
               & Single-$\gamma$ & 0.11844          & 0.854\%  & 0.10809          & 0.688\%  & 0.10316          & 0.900\%  & 0.10211          & 0.802\%  \\
               & \textit{Improve} & 14.72\%          &          & 12.94\%          &          & 13.78\%          &          & 14.86\%          &          \\
               & Multi-$\gamma$  & \textbf{0.11606} & 1.213\%  & \textbf{0.10582} & 0.709\%  & \textbf{0.10128} & 0.877\%  & \textbf{0.10094} & 0.804\%  \\
               & \textit{Improve} & 16.43\%          &          & 14.77\%          &          & 16.35\%          &          & 15.83\%          &          \\ \hline
5\%            & $\gamma=0$      & 0.18713          & 1.149\%  & 0.16040          & 0.844\%  & 0.14821          & 0.563\%  & 0.14313          & 0.356\%  \\
               & Single-$\gamma$ & 0.18594          & 0.590\%  & 0.15888          & 0.497\%  & \textbf{0.14495} & 0.588\%  & \textbf{0.13872} & 0.520\%  \\
               & \textit{Improve} & 0.64\%           &          & 0.95\%           &          & 2.20\%           &          & 3.08\%           &          \\
               & Multi-$\gamma$  & \textbf{0.18387} & 1.799\%  & \textbf{0.15834} & 1.028\%  & 0.14529          & 0.915\%  & 0.13933          & 0.738\%  \\
               & \textit{Improve} & 1.74\%           &          & 1.28\%           &          & 1.97\%           &          & 2.65\%           &          \\ \hline
\end{tabular}
\end{table}

\subsection{Reconstruction quality vs. depth.}
\begin{figure*}
\centering
\caption{The domain is divided into three regions with different depth.}
\label{fig:depth}
\includegraphics[width=0.5\textwidth]{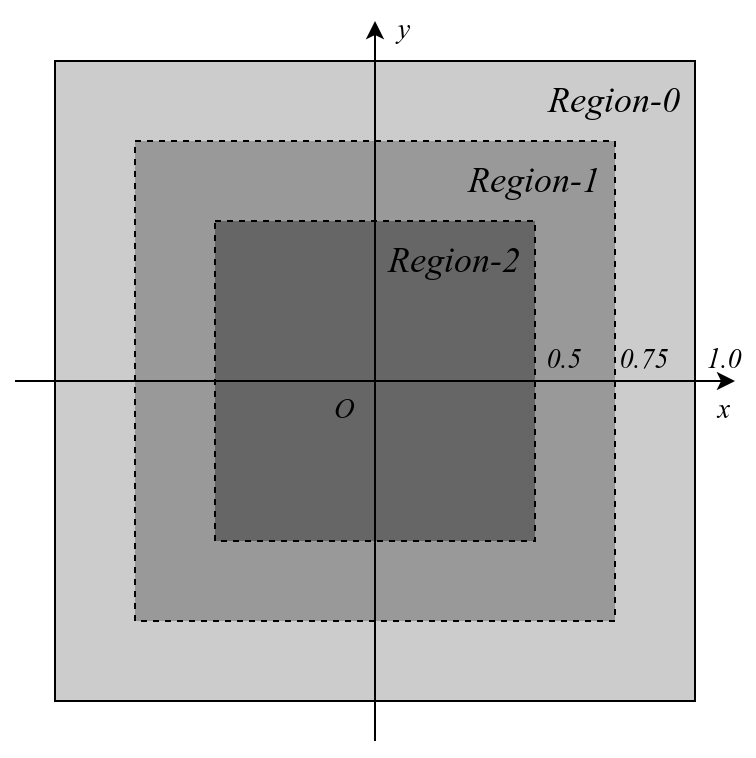}
\end{figure*}

From the experiments, the proposed method has demonstrated superior reconstruction quality (inclusion separation) from the boundary to the central area of the computational domain. These observations coincide with the traditional wisdom shown on depth-dependence estimates in classical inverse problem literature (e.g., see \cite{2009NagayasuUhlmannWangIP}).
Generally speaking, perturbations of the inclusions near the boundary have greater influence on boundary measurements, consequently any boundary data-based reconstruction model has the tendency to reconstruct inclusions near the boundary better.
To further study the performance of the proposed methods on regions with different depth,
we divide the domain $\Omega = [-1, 1]^2$ into three regions:
\begin{itemize}
  \item[(\textbf{0})] $[-1, 1]^2\setminus [-0.75, 0.75]^2$;
  \item[(\textbf{1})] $[-0.75, 0.75]^2\setminus [-0.5, 0.5]^2$;
  \item[(\textbf{2})] $[-0.5, 0.5]^2$.
\end{itemize}
As illustrated in Figure~\ref{fig:depth}, the depth grows from region \textbf{0} to region \textbf{2}.
The reconstruction results, generated by five random runs, on these three regions are shown in Table~\ref{tab:depth} (BCE) and Table~\ref{tab:depth_dice} (Dice),
showing that the learnable $\gamma$ can indeed improve the reconstruction on all different regions.

\begin{table}[htbp]
  \centering
  \caption{The test loss (binary cross entropy) on different regions of the
  domain with increasing depth.
  The reported mean and RSD are calculated based on 5 runs with different random seeds.
}
  \label{tab:depth}
\begin{tabular}{ll|ll|ll|ll}
\hline
\textbf{Noise} & \textbf{Model}  & \textbf{0}       &          & \textbf{1}       &          & \textbf{2}       &          \\
               &                 & mean             & RSD      & mean             & RSD      & mean             & RSD      \\ \hline
0\%            & $\gamma=0$      & 0.01117          & 1.137 \% & 0.05051          & 1.959 \% & 0.09818          & 2.146 \% \\
               & Single-$\gamma$ & \textbf{0.00764} & 3.600 \% & 0.03719          & 2.122 \% & 0.08153          & 2.895 \% \\
               & \textit{Improve} & 31.60\%          &          & 26.37\%          &          & 16.96\%          &          \\
               & Multi-$\gamma$  & 0.00769          & 2.643 \% & \textbf{0.03576} & 1.381 \% & \textbf{0.07939} & 3.508 \% \\
               & \textit{Improve} & 31.15\%          &          & 29.20\%          &          & 19.14\%          &          \\ \hline
1\%            & $\gamma=0$      & 0.01152          & 1.641 \% & 0.05480          & 1.821 \% & 0.10403          & 1.439 \% \\
               & Single-$\gamma$ & \textbf{0.00802} & 3.844 \% & 0.04243          & 1.295 \% & 0.09683          & 1.413 \% \\
               & \textit{Improve} & 30.38\%          &          & 22.57\%          &          & 6.92\%           &          \\
               & Multi-$\gamma$  & 0.00812          & 2.820 \% & \textbf{0.04174} & 1.019 \% & \textbf{0.09358} & 1.538 \% \\
               & \textit{Improve} & 29.51\%          &          & 23.83\%          &          & 10.05\%          &          \\ \hline
5\%            & $\gamma=0$      & 0.01370          & 1.561 \% & 0.06988          & 1.109 \% & \textbf{0.12325} & 1.874 \% \\
               & Single-$\gamma$ & \textbf{0.01165} & 0.904 \% & 0.06848          & 0.652 \% & 0.12684          & 2.140 \% \\
               & \textit{Improve} & 14.96\%          &          & 2.00\%           &          & -2.91\%          &          \\
               & Multi-$\gamma$  & 0.01257          & 1.749 \% & \textbf{0.06773} & 0.952 \% & 0.12462          & 2.053 \% \\
               & \textit{Improve} & 8.25\%           &          & 3.08\%           &          & -1.11\%          &          \\ \hline
\end{tabular}
\end{table}

\begin{table}[htbp]
  \centering
  \caption{The test loss (dice) on different regions of the
  domain with increasing depth.
  The reported mean and RSD are calculated based on 5 runs with different random seeds.
}
  \label{tab:depth_dice}
\begin{tabular}{ll|ll|ll|ll}
\hline
\textbf{Noise} & \textbf{Model}  & \textbf{0}       &          & \textbf{1}       &          & \textbf{2}       &          \\
               &                 & mean             & RSD      & mean             & RSD      & mean             & RSD      \\ \hline
0\%            & $\gamma=0$      & 0.20674          & 1.454 \% & 0.09955          & 1.197 \% & 0.17108          & 0.655 \% \\
               & Single-$\gamma$ & \textbf{0.16253} & 3.126 \% & 0.07500          & 2.327 \% & 0.14607          & 0.568 \% \\
               & \textit{Improve} & 21.38\%          &          & 24.66\%          &          & 14.62\%          &          \\
               & Multi-$\gamma$  & 0.16295          & 2.381 \% & \textbf{0.07406} & 2.003 \% & \textbf{0.14603} & 1.201 \% \\
               & \textit{Improve} & 21.18\%          &          & 25.61\%          &          & 14.64\%          &          \\ \hline
1\%            & $\gamma=0$      & 0.21417          & 1.084 \% & 0.11179          & 1.220 \% & 0.18801          & 0.661 \% \\
               & Single-$\gamma$ & 0.17080          & 2.817 \% & 0.08930          & 1.380 \% & 0.17990          & 0.731 \% \\
               & \textit{Improve} & 20.25\%          &          & 20.12\%          &          & 4.31\%           &          \\
               & Multi-$\gamma$  & \textbf{0.17069} & 2.394 \% & \textbf{0.08858} & 1.270 \% & \textbf{0.17568} & 0.934 \% \\
               & \textit{Improve} & 20.30\%          &          & 20.76\%          &          & 6.56\%           &          \\ \hline
5\%            & $\gamma=0$      & 0.24520          & 1.213 \% & 0.14875          & 0.951 \% & \textbf{0.23013} & 0.835 \% \\
               & Single-$\gamma$ & \textbf{0.21973} & 1.042 \% & 0.14575          & 0.627 \% & 0.23540          & 0.737 \% \\
               & \textit{Improve} & 10.39\%          &          & 2.02\%           &          & -2.29\%          &          \\
               & Multi-$\gamma$  & 0.22889          & 1.793 \% & \textbf{0.14472} & 0.774 \% & 0.23259          & 0.911 \% \\
               & \textit{Improve} & 6.65\%           &          & 2.71\%           &          & -1.07\%          &          \\ \hline
\end{tabular}
\end{table}

\subsection{Reconstruction on other datasets}

In this subsection, we present the numerical test of $\gamma$-deepDSM on more datasets with different conductivity distributions.
The first test concerns the inclusions with other shapes:
\begin{itemize}
\item [(\textbf{T1})] One triangle with three random vertices.
\item [(\textbf{T2})] Two equilateral triangles with random positions and scales.
\end{itemize}

The second test concerns continuous conductivity generated by multiple 2D
Gaussian distributions, rather than the previously-discussed piecewise-constant conductivity:
$$\begin{aligned}
  \sigma(\bm{x}) & =10^{J^*(\bm{x})}\in[1, 10], ~~~ J^*  :=\text{normalize}\left\{\sum_{k=0}^{K-1}J_k\right\}\in[0, 1], ~~K\in\mathbb{N}, \\
  J_k(\bm{x}) & =\exp\left\{-\frac{1}{2}(\bm{x}-\bm{c_k})^T\bm\Sigma_k^{-1}(\bm{x}-\bm{c}_k)\right\}, ~~~  \bm\Sigma_k^{-1}  =a_k^2 \begin{bmatrix}
    1-\mu_k^2\sin^2\theta_k & -\mu_k^2\sin\theta_k\cos\theta_k\\
    -\mu_k^2\sin\theta_k\cos\theta_k & 1-\mu_k^2\cos^2\theta_k
  \end{bmatrix},
\end{aligned}$$
where $\bm{c}_k, a_k, \mu_k, \theta_k$ are all sampled uniformly from
$[-0.9, 0.9]^2, ~~[2.0, 6.0], ~~[0.5, 1.0], ~~[0, 2\pi]$, respectively,
and the following two cases are considered:
\begin{itemize}
\item [(\textbf{G2})] Two Gaussian distributions, i.e., $K=2$.
\item [(\textbf{G3})] Three Gaussian distributions, i.e., $K=3$.
\end{itemize}
The reconstruction results are presented in Table \ref{tab:other_datasets} showing that learnable $\gamma$'s have
a similar enhancing effect.

\begin{table}[htbp]
  \centering
  \caption{The validation loss of the models trained on the data sets T1--G3.
The fLB operators can still be effective in inclusions with
  non-circular shapes and in continuous conductivity.}
  \label{tab:other_datasets}
\begin{tabular}{lllll}
\hline
\multicolumn{1}{c}{\multirow{2}{*}{\textbf{Model}}} & \multicolumn{4}{c}{\textbf{Datasets}}                                                                                                 \\ \cline{2-5} 
\multicolumn{1}{c}{}                                & \multicolumn{1}{c}{\textbf{T1}} & \multicolumn{1}{c}{\textbf{T2}} & \multicolumn{1}{c}{\textbf{G2}} & \multicolumn{1}{c}{\textbf{G3}} \\ \hline
$\gamma$=0                                          & 0.02717                         & 0.05447                         & 0.09534                         & 0.14294                         \\
Single-$\gamma$                                     & \textbf{0.02383}                & \textbf{0.04680}                & \textbf{0.06914}                & 0.12793                         \\
\textit{Improve}                                     & 12.29\%                         & 14.08\%                         & 27.48\%                         & 10.50\%                         \\
Multi-$\gamma$                                      & 0.02533                         & 0.04943                         & 0.07295                         & \textbf{0.12503}                \\
\textit{Improve}                                     & 6.77\%                          & 9.25\%                          & 23.48\%                         & 12.53\%                         \\ \hline
\end{tabular}
\end{table}

\subsection{Performance with varying number of data pairs}
\begin{figure*}[htbp]
  \centering
  \includegraphics[width=0.62\textwidth]{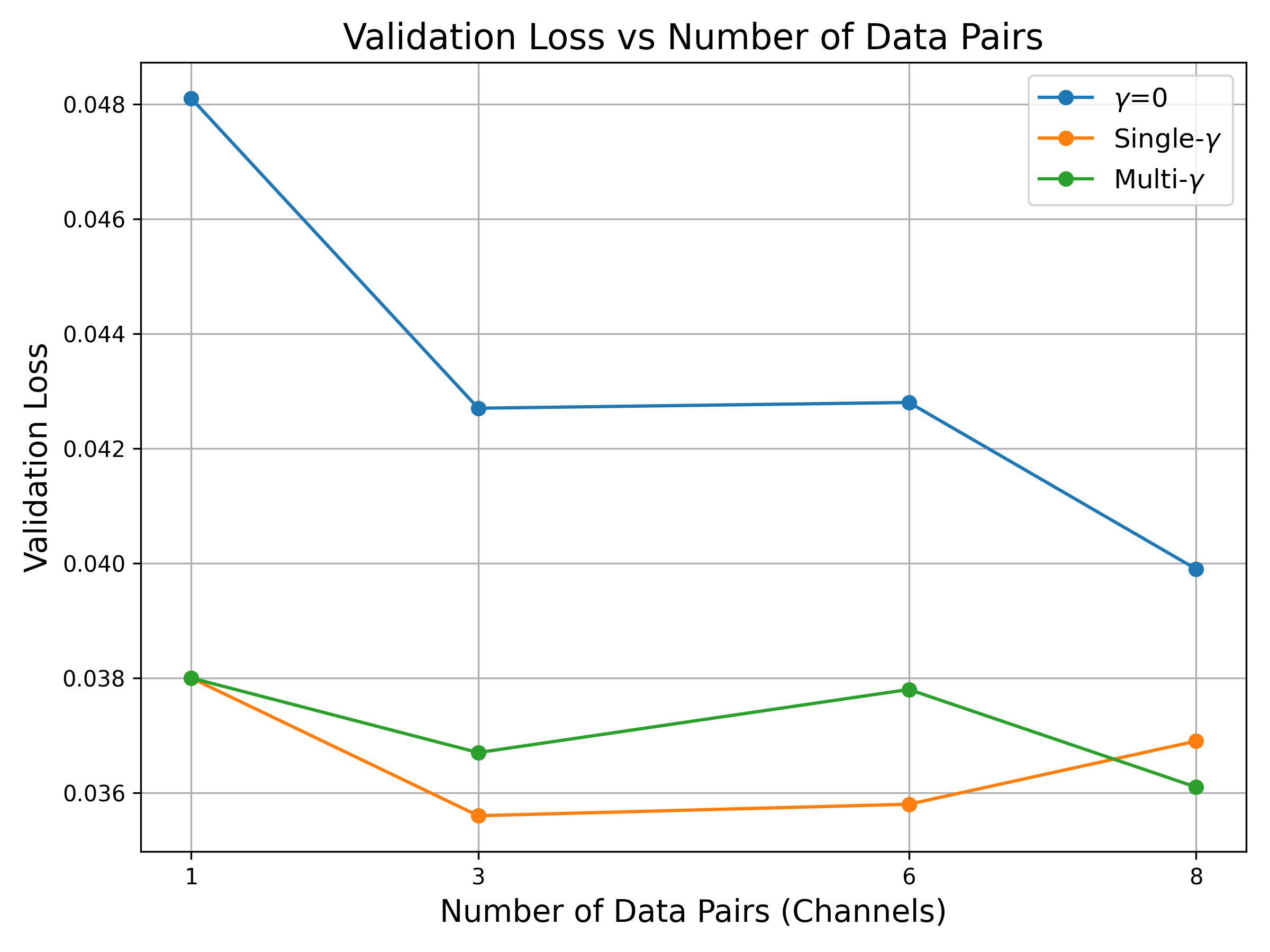}
  \caption{The validation results of the models with varying number of data pairs.
}
  \label{fig:data_pairs}
\end{figure*}

In this test, we evaluate the performance of $\gamma$-deepDSM with varying number of data pairs
by train the models with $L = 1, 3, 6$ data pairs, and evaluate the validation loss
(compared to the $L=8$ case).
Note that as the number of data pairs increases, the number of parameters in the
model also slightly increases. 
The results for varying number of data pairs are shown in Figure
\ref{fig:data_pairs}.
Adding more data pairs enriches information/features for a single sample (within expectation),
but it does not simply guarantee better performance,
as more data pairs with more learnable $\gamma$'s also increase the optimization complexity.
However, we also point it out that the $\gamma$-DDSM always performs bettern then the orginal DDSM,
even though the former one has fewer data pairs.

\subsection{Performance with lower spacial resolution}
With the same random seed, we also construct a lower-resolution but larger
version of the dataset \textbf{CIR3}: the conductivity distribution is defined
on a $32\times 32$ mesh instead of $64\times 64$, and the number of
samples is 4 times larger than the original one.
The results in Table~\ref{tab:low_res} show that the learnable $\gamma$'s can
still bring significant improvement.

\begin{table}[htbp]
  \centering
  \caption{The test loss on the lower-resolution dataset.
  The reported mean and RSD are calculated based on 5 runs with different random seeds.}
  \label{tab:low_res}
\begin{tabular}{ll|ll|ll}
\hline
$\delta$ & \textbf{Model}   & \textbf{BCE}     &          & \textbf{Dice}    &          \\
         &                  & mean             & RSD      & mean             & RSD      \\ \hline
0\%      & $\gamma=0$       & 0.03826          & 1.764 \% & 0.10775          & 2.169 \% \\
         & Single-$\gamma$  & 0.03135          & 0.752 \% & 0.08780          & 0.893 \% \\
         & \textit{Improve} & 18.06\%          &          & 18.52\%          &          \\
         & Multi-$\gamma$   & \textbf{0.03060} & 1.331 \% & \textbf{0.08695} & 0.930 \% \\
         & \textit{Improve} & 20.02\%          &          & 19.30\%          &          \\ \hline
1\%      & $\gamma=0$       & 0.04334          & 1.984 \% & 0.12342          & 1.793 \% \\
         & Single-$\gamma$  & 0.03872          & 1.661 \% & 0.11151          & 1.657 \% \\
         & \textit{Improve} & 10.66\%          &          & 9.65\%           &          \\
         & Multi-$\gamma$   & \textbf{0.03746} & 1.704 \% & \textbf{0.10841} & 1.826 \% \\
         & \textit{Improve} & 13.57\%          &          & 12.16\%          &          \\ \hline
5\%      & $\gamma=0$       & 0.05640          & 4.323 \% & 0.16223          & 3.180 \% \\
         & Single-$\gamma$  & \textbf{0.05571} & 1.375 \% & \textbf{0.16014} & 1.869 \% \\
         & \textit{Improve} & 1.22\%           &          & 1.29\%           &          \\
         & Multi-$\gamma$   & 0.05624          & 2.210 \% & 0.16276          & 2.417 \% \\
         & \textit{Improve} & 0.28\%           &          & -0.33\%          &          \\ \hline
\end{tabular}
\end{table}

\subsection{Comparison with other methods}

At last, we compare the performance of $\gamma$-deepDSM with a vanilla NN approach \cite{electronics7120422}, a PINN \cite{2025SmylTallmanHomaFlournoy},
and a Neural Inverse Operator (NIO) \cite{2023MolinaroYangEngquistMishra}.
We notice that the comparison against the original DDSM \cite{2020GuoJiang} ($\gamma=0$) has been conducted extensively above,
which clearly indicates the superior performance of learnable $\gamma$'s.

The vanilla NN approach maps from the boundary data
(voltage between adjacent electrodes) to the conductivity distribution.
It has 4 hidden layers with 150 neurons in each layer, resulting in nearly 0.98M
parameters, similar to the $\gamma$-deepDSM.
The leaky ReLU activation function is applied to each hidden layer.
We employ a learning rate $1.0$, weight decay of
$10^{-7}$. The model was trained for $100000$ epochs.

The PINN approach presented in \cite{2025SmylTallmanHomaFlournoy} contains 3
branches of fully-connected layers:
(1) for approximating the conductivity distribution,
(2) for approximating the boundary voltage (on electrode), and
(3) for approximating the voltage field in the domain.
To match the size of the $\gamma$-deepDSM, we set 22 neurons in each hidden layer, 
and the number of hidden layers for (1)-(3) to be 3, 2, 2, respectively.
The backbone of the these three branches consists of two convolutional blocks,
each with a 1D convolutional layer,
a batch normalization and a leaky ReLU activation function.
We employ an initial learning rate $10^{-1}$,
and a scheduler scaling the learning rate by $\frac{1}{\sqrt{\text{step}}}$,
warming up for 1000 steps.
The PINN was trained for 18600 epochs.

The NIO \cite{2023MolinaroYangEngquistMishra} is a novel architecture
for solving PDE inverse problems, which contains a DeepONet and
a Fourier Neural Operator (FNO).
We set its size to be 1.09M parameters, 
also close to the $\gamma$-deepDSM. 
We employ a initial learning rate of $10^{-3}$, and a
plateau learning rate scheduler with step size 15, decay factor 0.98.
The NIO was trained for 1000 epochs.

All models are trained and tested on a dataset of the same distribution as
\textbf{CIR3}. All data pairs (channels) were considered.
The results are presented in Table~\ref{tab:comparison}, showing that the $\gamma$-deepDSM with learnable
fLB order can achieve almost $\mathbf{50\%}$ improvement compared with the other models.

\begin{table}[htbp]
  \centering
  \caption{Binary cross entropy error and relative $L^2$ error evaluated on the
  testing samples.
  The relative $L^2$ error is defined as $\frac{\|\hat{\sigma}-\sigma^*\|_2}{\|\sigma^*\|_2}$,
  where $\hat{\sigma}$ is the reconstructed conductivity distribution, and $\sigma^*$ is the ground truth.}
  \label{tab:comparison}
\begin{tabular}{cccc}
\hline
                                 & \textbf{}       & \textbf{BCE $\downarrow$} & \textbf{Relative $L^2$ $\downarrow$} \\ \hline
\textbf{$\gamma$-deepDSM (ours)} & Single-$\gamma$ & 0.03553                   & 2.621\%                             \\
                                 & Multi-$\gamma$  & \textbf{0.03456}          & \textbf{2.557\%}                    \\
\textbf{NIO}                     &                 & 0.12070                   & 9.134\%                             \\
\textbf{ANN}                     &                 & 0.08350                   & 6.468\%                             \\
\textbf{PINN}                    &                 & 0.08246                   & 6.426\%                             \\ \hline
\end{tabular}
\end{table}

\section{Conclusion}

The present work contains two major contributions.
First, we have developed a ``solver-in-the-loop'' joint operator learning approach, $\gamma$-deepDSM, for solving the inverse problem of electrical impedance tomography, that is to reconstruct the media interfaces through the boundary measurements.
In this method, the ``solver-in-the-loop'' mechanism is adopted for incorporating PDE solvers with learnables for preconditioning.
Particularly, the scattered data on the boundary will first be processed through a learnable FEM module, and are lifted as high-dimensional features in the interior domain.
These features are used as different channels in the input of a U-Net, whose output is trained to reconstruct the conductivities in the computational domain.
The FEM module can learn the fractional orders of the Laplace-Beltrami applied to the boundary data through auto-differentiation, which significantly improves the reconstruction accuracy.
Second, to facilitate this process, we have implemented a tensorized FEM package, called LA-FEM, fully with PyTorch that enjoys the features of auto-differentiation, batch computation tradition of PyTorch to solve PDEs with multiple right-hand sides. The tensorized LA-FEM takes fully advantages of the CUDA parallelism by implementing the assembly and solving processes all in \texttt{einsum} using 4D tensors.
The proposed $\gamma$-deepDSM enabled by LA-FEM capitalizes on the strengths of both traditional FEM (highly accurate) and operator learning (capturing the nonlinear relationship from distribution level), mitigating their respective limitations. This synergistic approach expands new possibilities to tackle difficult inverse problems, instead of trying to replace traditional numerical solvers with neural networks.
The source code of the present work including LA-FEM is available at \url{https://github.com/weihuayi/fealpy/tree/master/app/lafem-eit}.

\section{Appendix}
\label{appendix}

In this section, we shall present the development details of LA-FEM and some guidelines for using this package.

\subsection{FEM based on Tensor computing}

Our key technique for implementing FEM in PyTorch is to reformulate all its computational components in a fully tensorized manner.
For presentation purpose, we resort to the Laplace equation to illustrate how this can be done.
Given a computation domain $\Omega$ with arbitrary dimension $D$,
we consider
\begin{equation}
\label{lap_eq1}
(\nabla u, \nabla v)_{\Omega} = (f,v) ~~~~~ \forall v\in H^1_0(\Omega), ~~~~ u = 0 ~~~~~ \text{on} ~ \partial\Omega.
\end{equation}

Here we consider the FEM of the primal form of \eqref{lap_eq1}. 
The key idea of FEM is to generate a mesh of the computation domain $\Omega$ and a set of basis functions associated with this mesh $\{ v_1, ..v_R \}$.
Then, the associated discretization scheme reads
\begin{equation}
\label{lap_eq2}
(\nabla u_h, \nabla v_r)_{\Omega} = (f,v_r)_{\Omega} ~~~~~ \forall v_r ,
\end{equation}
where the zero boundary condition is immersed in the functions $v_r$.
As the elliptic equations are considered, we use the conforming Lagrange elements for those basis functions
in which each $v_r$ enjoys a local support property, i.e., its support usually contains the elements surrounding one mesh node.
Denote the Lagrange FE space as $\mathbb{V}_h= \textrm{Span}\{ v_1, ..v_R \}$.
For more details, we refer readers to \cite{2008BrennerScott}.
The implementation of FEM can be divided into several key modules. 
Here, we focus our tensor computing strategy on two critical tasks: 
generating basis functions and assembling stiffness matrices.
We also refer readers to \cite{2008Chen} for an FEM package with its notes for an introduction to tensor computing of FEM.

Tensorizing generation of basis functions heavily relies on the barycentric-coordinate representation of Lagrange shape functions \cite{2008BrennerScott,2014RandGilletteBajaj}.
For the convenience of the readers, this coordinate system is briefly reviewed here.
Consider a $D$-dimensional simplicial mesh $\mathcal{T} = [T_1, \cdots, T_N]$ of $\Omega$.
Each element $T_n$ in $\mcT$ is a $(D+1)$-simplex:
$$
T_n = \text{ConvexHull}(\bfv_{n0}, \bfx_{n1}, \dots, \bfx_{nD}) := \left\{\sum_{b=0}^{D}\lambda_{b}\bfx_{nb}: 0\le \lambda_{b} \le 1, \sum_{b=0}^{D}\lambda_{b} = 1\right\},
$$
where $\{\bfx _{nb}\}_{b=0}^D$ denotes the vertices of $T_n$.
For any $\bfx \in T_n$, its unique barycentric coordinates $\blambda_n :=(\lambda_{n0}, \lambda_{n1}, \cdots, \lambda_{nD})$ satisfy
$$
\bfx = \sum_{b=0}^D \lambda_{nb}\bfx_{nb}, ~~~~ \sum_{b=0}^D \lambda_{b} = 1, \text{ and } \lambda_{b} \ge 0.
$$
As such, $\lambda_{nb}(\cdot)$ can be viewed a mapping from $\bfx$ to $\lambda_{nb}$.

On each element $T_n$, we denote the set of Lagrange shape functions by $\{\phi_{nk}\}_{k=1}^K$,
where $K$ depends on the dimension $D$ and the polynomial degree.
One crucial property beneficial to tensor computing is that the shape functions admit an expression of barycentric coordinates: 
$$
\phi_{nk}(\bfx) = f_k(\blambda_{n}(\bfx)),
$$
where each $f_k(\lambda_0,...,\lambda_D)$ is a $(D+1)$-variable function that is fixed once the finite element space is chosen.
Then, we select $Q$ quadrature points for every element that share the same barycentric coordinates $\blambda_{q}$ with the weights $\{w_{q}\}$, $q=1,\dots,Q$. We also note that their corresponding Cartesian coordinates $\{ \bfx_{nq} \}$ depend on $T_n$, which are different from elements to elements.

\begin{figure}[h]
    \centering
    \includegraphics[width=0.8 \textwidth]{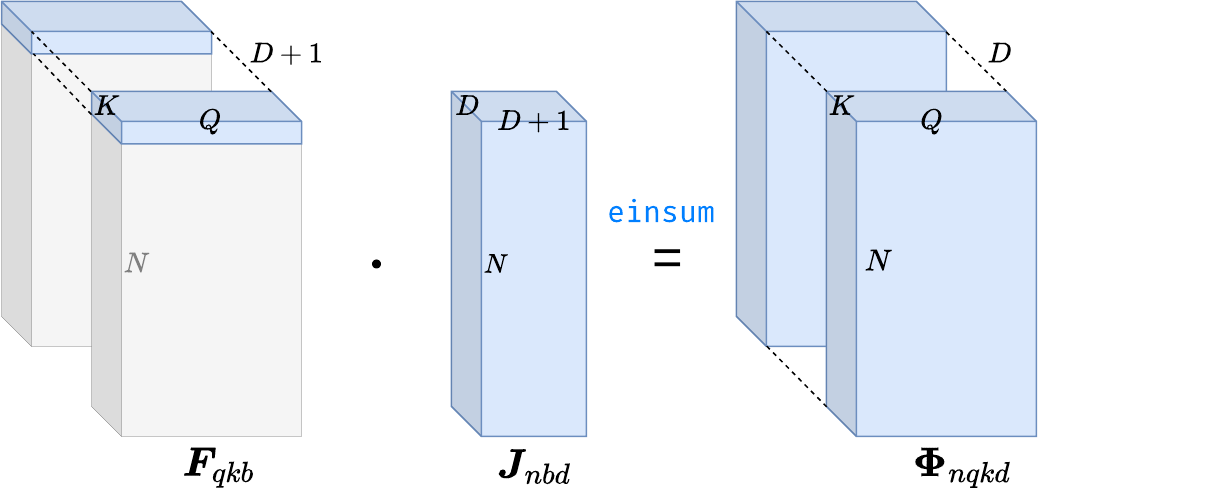}
    \caption{The flowchart of how the Mesh module generates gradient of shape functions.}
    \label{fig:tensor_grad_phi}
\end{figure}

We first represent the gradient of shape functions by a 4D tensor $\bPhi_{nqkd}$, $1\le n \le N$, $1\le q \le Q$, $1 \le k \le K$, $1 \le d \le D$.
With the \textbf{Einstein summation} convention, 
this can be done via the chain rule: 
\begin{equation}
\label{tensor_grad}
\bPhi_{nqkd} := \partial_{x_d} \phi_{nk}(\bfx_{nq})
=\frac{\partial f_{k}}{\partial \lambda_{b}}(\blambda_{q})
\frac{\partial \lambda_{nb}}{\partial x_d},
\end{equation}
where the summation is taken along the $b$-axis that is repeated.
From \eqref{tensor_grad}, we denote $\bfF_{qkb} := \frac{\partial f_{k}}{\partial \lambda_{b}}(\blambda_{q})$ 
and a batched Jacobian tensor $\bfJ_{nbd} := \frac{\partial \lambda_{nb}}{\partial x_d} $, which have 
the shapes of $[Q, K, D+1]$ and $[N, D+1, D]$, respectively.
Notice that the values of $f_{k}$ and its derivatives at $\blambda_{q}$ are independent of elements $T_n$.
As a result, $\bfF_{qkb}$ is independent of the geometry because all elements share the same barycentric coordinates for quadrature points.
In FEALPy, the \texttt{Mesh} module can evaluate the gradient of shape functions with respect to the barycentric coordinates, i.e., $\bfF_{qkb}$, on the reference element. 
$\bfF_{qkb}$ is independent of the shape of the element, which only needs to be constructed for the reference element.  
Then, $\bfF_{qkb}$ can be broadcasted in to a tensor for all the elements in the mesh with the shape of $[N, Q, K, D+1]$ in subsequent tensor calculations through \eqref{tensor_grad}. 
See Fig. \ref{fig:tensor_grad_phi} for an illustration.

\begin{figure}[h]
    \centering
    \includegraphics[width=0.8\textwidth]{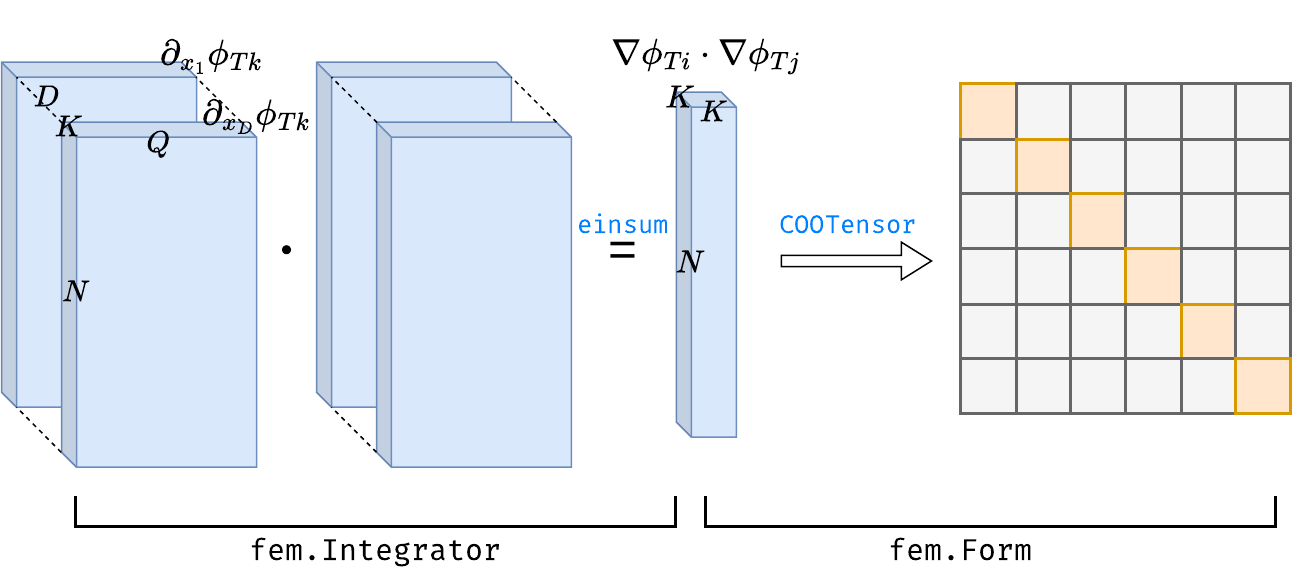}
    \caption{Local integration to sparse matrix assembly. 
    The weights and element measures are omitted in this illustrative diagram for brevity.}
    \label{fig:tensor_mat}
\end{figure}

The next stage is matrix assembly, and we need to compute
\begin{equation}
\label{intgrand}
\sum_{d=1}^D\int_{T_n} \partial_{x_d} \phi_{nk} \partial_{x_d} \phi_{nk} \dd x 
\approx \sum_{d=1}^D \sum_{q=1}^Q |T_n| w_q \partial_{x_d} \phi_{nk}(\bfx_{nq})  \partial_{x_d} \phi_{nk}(\bfx_{nq}).
\end{equation}
In terms of tensors, computing \eqref{intgrand} involves the contraction of tensor $\bPhi_{nqkd}$ and itself to output a tensor shaped $[N, K, K]$,
where the summation is taken along both the $d$-axis (physical dimension), and a $w_q$-weighted inner product is taken along the $q$-axis (number of quadrature points).
This procedure is implemented in the \texttt{Integrator} class from the \texttt{FEM} module.
By the mapping from local degrees of freedom (DoFs) to global DoFs provided by the \texttt{FunctionSpace} module, the \texttt{Form} class from the \texttt{FEM} module manages to assembly the global sparse matrix. 
See Fig. \ref{fig:tensor_mat} for the illustration.
In this process, we rely on a fast and optimized implementation of \texttt{einsum} in PyTorch, as well as other backends, to be the major tool used in the \texttt{Integrator} class. \texttt{einsum} offers convenient access to a wide range of tensor operations in a compact and readable manner,
particularly for higher-dimensional tensors and multiple inputs.
The \texttt{COOTensor} class is our FEALPy implementation of the coordinate sparse format, which supports multiple tensor computing backends, and it will be discussed in the next subsection.
In matrix assembly, they will be converted to the Compressed Sparse Row (CSR) format by default, which is more commonly used in numerical linear algebra packages such as numpy.

We employ the techniques discussed above to update FEALPy on multiple tensor computing backends.
Our goal is to design such an FEM module that incorporates the three important features: AD, tensorized computation taking full advantage of CUDA, and batched computation of PDEs. 
The traditional way of looping through elements/vertices with conditionals on indices will break the AD. A native and backend-agnostic implementation (e.g., the Einstein sum is called as \texttt{bm.einsum}) in support multiple backends help avoid the data format conversions (e.g., between \texttt{np.ndarray} and \texttt{torch.Tensor}, the procedure of which needs manually pass the AD's values).
Most importantly, the tensorized assembly avoids the computation of geometry-hard-coded derivatives for different finite element shape functions, which yet again will break AD.
LA-FEM is implemented in FEALPy with these potential behavior-broken software designs in mind. 
We defer further details at the software design level to the next subsection.

\subsection{Software Design of LA-FEM in FEALPy}
\label{sec:lafem_impl}

\subsubsection{Backend Manager for Tensor Computing}
\label{appendix:backend-manager}
One of the most important interface in FEALPy is a backend system lying in the Tensor Level shown in Fig. \ref{fig:fealpy_structure}.  Inspired by Tensorly
\cite{Jean2019tensorly} that has a backend system supporting multiple backends
such as NumPy, PyTorch, TensorFlow, and others, we design a similar backend
system that incorporates nearly all standard Application Programming Interfaces
(APIs) to fulfill the needs from all the high-level modules in FEALPy.

\begin{figure}[htbp]
    \centering
    \includegraphics[width=0.8\textwidth]{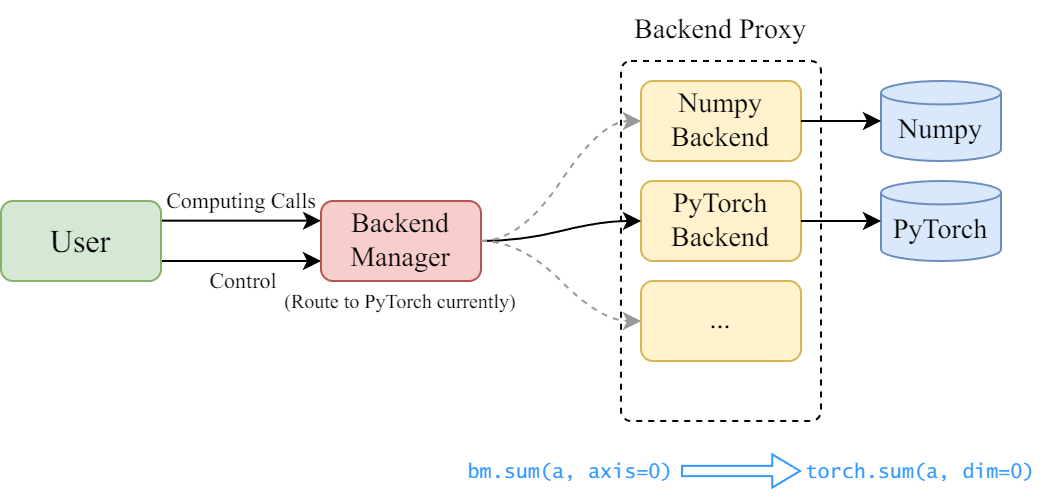}
    \caption{Design of the backend system in FEALPy. The backend manager object is a ``router'', and the backend proxy objects are ``translators''.}
    \label{fig:backend}
\end{figure}

The core of FEALPy's backend system is an object called \textbf{Backend Manager} (\texttt{BM}), which
centrally mageans various computing backends and serves as the sole interface
for users to interact with the Tensor Level.  All higher-level modules operate
solely on \texttt{BM} to access backend APIs.  Below is a simple Python code to
import the Backend Manager:
\begin{lstlisting}
from fealpy.backend import backend_manager as bm
\end{lstlisting}
\texttt{BM} defines a series of tensor calculation functions mainly based on the ``Python array API standard v2023.12'' \cite{PythonArrayAPI}.
For instance, the function \verb|bm.zeros| generates a tensor of a specified shape filled with zeros.
\texttt{BM} serves as a wrapper around several computing backends (e.g., NumPy and PyTorch) with little to no overhead, and has the capability to switch between them.
The default computing backend is NumPy, meaning that the \verb|bm.zeros| instruction refers to \verb|np.zeros|, 
and the returned object is a NumPy array.
The computing backend can be switched by the \verb|set_backend| function, for example the code below
\begin{lstlisting}
bm.set_backend('pytorch')
\end{lstlisting}
set the default computing backend as PyTorch. 
Then, \verb|bm.zeros| refers to \verb|torch.zeros|, returning a PyTorch tensor,
and all higher-level modules operate solely on PyTorch, leveraging its
tensor computing capabilities.
To achieve this user-friendly function,
we have designed a special backend proxy communicating with each backend.
When the user calls the tensor computing API to \texttt{BM}, 
it invokes the corresponding API in the selected backend proxy, 
translating it to the specific computing backend, as illustrated in Fig. \ref{fig:backend}.

\subsubsection{Batched computation}
In the context of problem \eqref{AD_eq1} in the proposed method, it requires solving $L$ equations. These equations correspond to solving a boundary value problem $\Gamma = \{\gamma_1, ..., \gamma_L \}$ with $L$ different boundary conditions, respectively. This process necessitates solving PDEs with
multiple right-hand sides vectors $\bfA u_h = \bfF,$ where the $\bfF$ is a
matrix with $L$ columns. The software design to cater this behavior is similar
to those in batched computations during NN training, meaning that the package
must be capable of solving multiple PDEs simultaneously in a parallelized
manner. 
To support the batched computation, the \texttt{Integrator} objects in FEALPy are specifically designed to accommodate coefficients and sources with a batch
dimension along the first dimension ($0$-axis) conforming with PyTorch/JAX
tradition.  Users can specify the \verb|batched| and \verb|batch_size|
parameters for \texttt{Integrator} and \texttt{Form} respectively, as shown in
the following code:
\begin{lstlisting}
    bform = BilinearForm(space, batch_size=50)
    bform.add_integrator(ScalarDiffusionIntegrator(coef=coef, batched=True))
\end{lstlisting}
This feature ensures that the output tensor retains the batch dimension yield
the associated FEM solutions.

\subsubsection{Auto-differentiation}
In addition to Tensor computing, one feature we need for DL is auto
differentiation.  The auto differentiation system of some backends, like PyTorch
and JAX, helps to calculate the gradient and even higher-order derivatives
automatically with the help of computational graph.
Our goal in this context is to be able to calculate both spatial gradient and
gradient of a parameter. Specifically, $\partial_{\Gamma} \Phi$ in \eqref{AD_eq1} for any $\Gamma$, and $\partial_{\gamma}\mathcal{L}$ in \eqref{optim_problem} are needed.
The proposed implementation gives users a generic interface, which is not
necessarily only for the fractional order studied for the EIT problem. This
implementation enables gradient computations for sources, coefficients of PDEs,
mesh nodes coordinates, and even have the capacity to parametrize shape
functions.  The critical challenge is to maintain certain topologies and the
differentiability of the computational graph generated by PyTorch, in which
nodes represent either tensors or an operation.  All user-implemented
operations/functions will be automatically registered in this computational
graph. In our method, these operations/functions, upon being linked with
subsequent NN modules, must retain the differentiability with respect to the
parameters that reside in the FEM solutions implicitly.  This requirement
demands significant differences in our code compared to typical NumPy
implementations, since discontinuous and non-differentiable operations is ought
to be avoided, for example, sorting or permuting the parameters or coordinates
whose gradients need to be tracked, retrieving indices from the maximum element
in an iterable, in-place operations that manually set values, all of which are
quite common in FEM assembly procedure.

In summary, all these three features significantly enhance the compatibility of
FEM with DL in terms of data structure, forward propagation, and back propagation.
To fully utilize the strengths of DL-oriented backends such as PyTorch or Jax, the native implementation eliminates the needs for data conversion, thus facilitates more much efficient data I/O with Neural Network models.

\subsubsection{Sample tests for LA-FEM}

In our test of the multi-backend mat-mat multiplication for LA-FEM,
the backend is specified through the function \verb|bm.set_backend|,
and there is no need to modify any other code after the backend selection.
See Listing \ref {lst:matmul}.
\begin{lstlisting}[label=lst:matmul, 
caption={Backend selection for tensor computing.}]
  from fealpy.backend import backend_manager as bm

  bm.set_backend('numpy')

  if bm.backend_name == 'pytorch':
      bm.set_default_device('cuda')

  A = bm.random.rand(SIZE, SIZE)
  F = bm.random.rand(SIZE, SIZE)
  A = bm.astype(A, bm.float64)
  F = bm.astype(F, bm.float64)
  u = bm.matmul(A, F)
\end{lstlisting}

In the test of the batched FEM for LA-FEM, \textbf{Integrator} and \textbf{Form}
objects should accept multiple boundary conditions and cell matrix, respectively.
List \ref{lst:batch_test} is an example for a batched right-hand side assembly in a spatial dimension- and backend-agnostic way:
\begin{lstlisting}[label=lst:batch_test,
    caption={Multiple PDEs can be solved in parallel, different PDEs are treated as different instances in a batch. After the assembly on a given mesh or difference meshes with the same dimensions, the user only needs to change the batch dimension when incorporating new source or boundary conditions.}]
  lform = LinearForm(space, batch_size=10)
  lform.add_integrator(ScalarSourceIntegrator(pde.source, batched=True))
  lform.add_integrator(ScalarNeumannBCIntegrator(pde.neumann, batched=True))
  F = lform.assembly()
\end{lstlisting}

\end{document}